\documentclass[a4paper, 10pt]{amsart}
\usepackage{amsthm}
\usepackage[]{amsmath}
\usepackage{amssymb}
\usepackage{enumerate}
\usepackage{tabularx}
\usepackage[]{color}
\usepackage[left=2.8cm, right=2.8cm, top=2.5cm, bottom=2.5cm]{geometry}
\usepackage[colorlinks]{hyperref}
\usepackage{tikz}
\usepackage{multirow}
\usepackage{diagbox}
\usepackage{subcaption}
\usepackage[ruled, vlined]{algorithm2e}
\usepackage{xcolor}

\allowdisplaybreaks[4]

\title[Preconditioned FEM with $C^0$RA]{A Preconditioned Discontinuous
Galerkin Method for Biharmonic Equation with $C^0$-Reconstructed
Approximation}

\author[R. Li]{Ruo Li} \address{CAPT, LMAM and School of Mathematical
Sciences, Peking University, Beijing 100871, P.R. China}
\email{rli@math.pku.edu.cn}

\author[Q.-C. Liu]{Qicheng Liu} \address{School of Mathematical
Sciences, Peking University, Beijing 100871, P.R. China}
\email{qcliu@pku.edu.cn}

\author[F.-Y. Yang]{Fanyi Yang} \address{College of Mathematics
, Sichuan University, Chengdu 610065, P.R. China}
\email{yangfanyi@scu.edu.cn}

\newcommand{\bm}[1]{\boldsymbol{#1}}
\newcommand{\bmr}[1]{\bm{\mr{#1}}}
\newcommand{\lj}{[ \hspace{-2pt} [}
\newcommand{\rj}{] \hspace{-2pt} ]}
\newcommand{\mb}[1]{\mathbb{#1}}
\newcommand{\mc}[1]{\mathcal{#1}}
\newcommand{\mr}[1]{\mathrm{#1}}
\newcommand{\jump}[1]{\lj #1 \rj}
\newcommand{\aver}[1]{ \{#1\}  }
\newcommand{\wt}[1]{ \widetilde{ #1}}
\newcommand{\wh}[1]{ \widehat{ #1}}

\newcommand{\DGenorm}[1]{ \| #1\|_{{e}}}

\newcommand{\DGtenorm}[1]{|\!|\!| #1 |\!|\!|_{e}}
\newcommand{\snorm}[1]{ \| #1 \|_{s}}

\newcommand{\dx}[1]{\mathrm{d} \boldsymbol{#1}}

\def\dim{\ifmmode \mathrm{dim} \else \text{dim}\fi}
\def\rank{\ifmmode \mathrm{rank} \else \text{rank}\fi}
\def\MTh{\mc{T}_h}

\def\MEh{\mc{E}_h}
\def\MNh{\mc{N}_h}

\def\MEhB{\mc{E}_h^{\backslash \Gamma}}
\def\un{\bm{\mr n}}
\def\lap{\Delta}

\def\mI{\mc{I}}
\def\mS{\mc{S}}
\def\mN{\mc{N}}
\def\mD{\mc{D}}
\def\mR{\mc{R}}
\def\mP{\mc{P}}

\def\mV{\mc{V}}

\def\MEhI{\mc{E}_h^i}
\def\MEhB{\mc{E}_h^b}

\def\mE{\mc{E}}

\def\MNh{\mc{N}_h}
\def\MNhI{\mc{N}_h^i}
\def\MNhB{\mc{N}_h^b}

\newtheorem{assumption}{Assumption}
\newtheorem{theorem}{Theorem}
\newtheorem{lemma}{Lemma}

\newtheorem{remark}{Remark}

\definecolor{orange}{rgb}{1, 0.5, 0}

\begin{document}

\maketitle

\begin{abstract}
In this paper, we present a high-order finite element method based on
a reconstructed approximation to the biharmonic equation.
In our construction, the space is reconstructed from nodal values 
by solving a local least squares fitting problem per element.
It is shown that the space can achieve an arbitrarily high-order
accuracy and share the same nodal degrees of freedom with the $C^0$
linear space.
The interior penalty discontinuous Galerkin scheme can be directly
applied to the reconstructed space for solving the biharmonic
equation.
We prove that the numerical solution converges with optimal orders
under error measurements.
More importantly, we establish a norm equivalence between the
reconstructed space and the continuous linear space. 
This property allows us to precondition the linear system arising from
the high-order scheme by the continuous linear space on the same mesh. 
This preconditioner is shown to be optimal in the sense that the
condition number of the preconditioned system admits a uniform upper
bound independent of the mesh size.
Numerical examples in two and three dimensions are provided to
illustrate the accuracy of the scheme and the efficiency of the
preconditioning method.

\noindent \textbf{keywords}: 
reconstructed approximation;
least squares fitting;
interior penalty discontinuous Galerkin method;
low-order preconditioning; 
multigrid method;

\end{abstract}


\section{Introduction}
The biharmonic equation is a fourth-order elliptic equation arising in
fields of continuum mechanics.
It models the thin plate bending
problem in continuum mechanics and describes the slow flows of
viscous incompressible fluids. 
For numerically solving the biharmonic equation, 
finite element method has been a standard numerical technique, 
and till now, there are many successful finite element methods for
this problem, known as conforming, non-conforming and mixed methods.

Conforming methods use $C^1$ elements in the scheme, which satisfy the
$C^1$ continuity conditions across interelement boundaries; see
\cite{Argyris1968tuba, Douglas1979family, Ciarlet2002finite}. 
The implementation of such elements and spaces is quite complicated
even in two dimensions. 
The non-conforming methods relax the strong continuity requirements
across the interelement faces and avoid the construction of $C^1$
elements; we refer to \cite{Adini1960analysis, Morley1968triangular,
Rannacher1979nonconforming, Shi1986convergence} for some typical
elements.
The construction of such elements with high-order polynomials is
also non-trivial, and they are rarely used in practice.
The mixed finite element method (see \cite{Boffi2013mixed} and the
references therein) is another standard method in
solving the biharmonic equation, which also avoids using $C^1$ 
elements.
The fourth-order equation is rewritten into a system of second-order
elliptic equations, and the mixed system can be solved with $C^0$
spaces. 
The mixed method introduces an auxiliary variable, which increases
the number of unknowns and leads to a saddle point linear system. 
In recent years,
the discontinuous Galerkin (DG) finite element methods and
the $C^0$ interior penalty methods have been developed for the
biharmonic equation; see \cite{Mozolevski2003priori,
Suli2007hp, Mozolevski2007hp, Georgoulis2008discontinuous,
Cockburn2009hybridizable, Gudi2008mixed} for DG methods and see
\cite{Engel2002continuous, Brenner2005interior, Brenner2015C0} for
$C^0$ interior penalty methods.
Both methods can also be regarded as non-conforming methods. 
The $C^1$ continuity conditions are weakly imposed by applying proper
penalties on the interelement faces in the bilinear form. 
Consequently, simple finite element spaces such as discontinuous
piecewise polynomial spaces and $C^0$ Lagrange finite element spaces
can be used to solve the biharmonic equation. 
The elements in both spaces are much simpler than $C^1$ elements
especially when using high-order polynomials.
In addition, the DG methods allow totally discontinuous elements, 
which provide great flexibility in the mesh partition but introduce
more degrees of freedom than $C^0$ spaces. 
To achieve a high-order accuracy, an alternative approach is to apply
some reconstruction techniques in the scheme, which may give a better
approximation but with less degrees of freedom. 
In \cite{Lamichhane2014finite, Chen2017C0, Guo2018C0,
Huang2020recovery}, the recovery based methods were developed and
analysed, which solved the biharmonic equation by $C^0$ linear
elements.
The main idea of such methods is to employ a gradient recovery
operator to reconstruct the piecewise constant gradient as a piecewise
linear function.
Then, computing the second-order derivatives for $C^0$ linear finite
element functions becomes possible.
In \cite{Li2017discontinuous, Li2016discontinuous}, the authors
proposed a reconstructed finite element method to the biharmonic
equation, where the high-order approximation space was obtained by a
patch reconstruction process from the piecewise constant space. 
The idea of reconstructing high-order spaces from low-order spaces
in solving the fourth-order equation can be traced back to 
the works on developing plate and shell
elements; see \cite{Onate1993derivation, Phaal1992simple}.

In this work, we extend the idea of the reconstructed method in
\cite{Li2016discontinuous} to the biharmonic equation. 
The high-order approximation space is reconstructed from the $C^0$
linear space, i.e.  only nodal degrees of freedom are involved. 
For each element, we construct a wide patch and solve a local least
squares fitting problem from nodal degrees of freedom in the patch to
seek a high-order polynomial. 
We prove that any high-order accuracy can be achieved by the
reconstructed space while the space shares the same degrees of freedom
with the continuous linear space.
The new space consists of discontinuous piecewise polynomial
functions, and can be regraded as a small subspace of the standard
discontinuous piecewise polynomial space. 
Thus, it is natural to employ the interior penalty discontinuous
Galerkin scheme to the reconstructed space for solving the biharmonic
equation, and the error estimation is also straightforward under the
Lax-Milgram framework.
We prove the optimal convergence rates under error measurements, which
are also confirmed by a series of numerical tests in two and three
dimensions. 

For the fourth-order equation, solving the resulting linear system
efficiently is of particular concern in finite element methods because
the final matrix is extremely ill-conditioned, which, generally, has a
condition number of $O(h^{-4})$.
Therefore, it is very desired to develop an efficient linear solver. 
The most common methods are to construct a proper preconditioner,
which can remarkably reduce the condition number in the
preconditioned system, such as multigrid methods and domain
decomposition methods. 
For classical conforming and non-conforming elements, 
the preconditioning methods have been well developed; see
\cite{Zhang1996two, Zhang1989optimal, Brenner1996two,
Brenner1989optimal, Tang2017local, Stevenson2003analysis,
Carstensen2021hierarchical}.
As mentioned in \cite{Karakashian2018two}, there are few works
concerning the preconditioning techniques for the penalty methods. 
In \cite{Karakashian2018two, Feng2005two}, the authors presented
two-level non-overlapping and overlapping additive Schwarz methods for
DG methods, where the bound of the condition number to the
preconditioned system was given. 
In \cite{Smears2018nonoverlapping, Smears2014discontinuous}, the
authors developed the non-overlapping domain decomposition method for
DG methods solving HJB equations. 
For $C^0$ interior penalty methods, we refer to \cite{Brenner2005two}
and \cite{Brenner2006multigrid} for the two-level Schwarz
method and the multigrid method, respectively.
For the proposed method, the main feature is that
the reconstructed space uses only nodal degrees of freedom to achieve
a high-order accuracy.
From this fact, we can show that the inverse of the matrix arising
from the lowest-order scheme, i.e. the scheme over the continuous
linear spaces, can serve as a preconditioner to any high-order
reconstructed space on the same mesh.
The low-order preconditioning is also a classical technique in finite
element methods for preconditioning the high-order scheme; see
\cite{Chalmers2018low, Pazner2023low, Li2023preconditioned} for some
examples. 
We establish a norm equivalence between the reconstructed
space and the continuous linear space. 
This crucial property allows us to prove that the preconditioner is
optimal by showing the condition number of the preconditioned
system admits a upper bound independent of the mesh size.
Moreover, we propose a multigrid method as an approximation to the
inverse of the matrix from the lowest-order scheme, following the idea
in \cite{Vanek2001convergence}. 
The convergence analysis to the multigrid method is also presented.
The numerical tests in two and three dimensions illustrate the
efficiency of the preconditioning method.

The rest of this paper is organized as follows. 
In Section \ref{sec_preliminaries}, we give some notation used in the
scheme.
In Section \ref{sec_space}, 
we introduce the reconstructed approximation space and prove the
basic properties of the space.
Section \ref{sec_scheme} presents the interior penalty scheme for the
biharmonic equation. 
The preconditioning method is also proposed and analysed in this
section.
In Section \ref{sec_numericalresults}, we conduct a series of
numerical tests to demonstrate the accuracy to the proposed scheme
and the efficiency to the preconditioning method.


\section{Preliminaries}
\label{sec_preliminaries}
Let $\Omega \subset \mathbb{R}^d (d = 2, 3)$ be a bounded convex
polygonal (polyhedral) domain with the boundary $\partial \Omega$.
Let $\MTh$ be a quasi-uniform triangulation of the domain $\Omega$
into triangular (tetrahedral) elements.
For any $K \in \MTh$, we let $h_K := \text{diam}(K)$ be its
diameter, and let $\rho_K$ be the radius of the largest ball
inscribed in $K$. 
We let $h := \max_{K \in \MTh} h_K$ and $\rho := \min_{K \in \MTh}
\rho_K$, where $h$ is also the mesh size.
The mesh $\MTh$ is assumed to be quasi-uniform in the sense that
there exists a constant $C_\nu$ independent of $h$ such that $h \leq
C_\nu \rho$, where $\rho := \min_{K \in \MTh} \rho_K$.
For any $K \in \MTh$, we let $w(K) := \{K' \in \MTh: \ \overline{K}
\cap \overline{K'} \neq \varnothing \}$ be the set of elements
touching $K$.
Let $\MEh$ be the set of all $d-1$ dimensional faces in $\MTh$, and we
decompose $\MEh$ into $\MEh = \MEhI + \MEhB$, where $\MEhI$ and
$\MEhB$ are collections of interior faces and the faces lying on the
boundary $\partial \Omega$, respectively.
For any $e \in \MEh$, we let $h_e := \text{diam}(e)$ be its diameter. 
We further denote by $\MNh$ be the set of all nodes in 
$\MTh$. 
Similarly, $\MNh$ is decomposed into $\MNh = \MNhI + \MNhB$, where
$\MNhI := \{ \bm{\nu} \in \MNh: \ \bm{\nu} \in \Omega\}$ and $\MNhB :=
\{ \bm{\nu} \in \MNh: \ \bm{\nu} \in \partial \Omega \}$. 
For any $K \in \MTh$, we define  $\mc{N}_K := \{ \bm{\nu} \in \MNh:
\bm{\nu} \in \partial K\}$ as the set of all vertices of the element
$K$.

Next, we introduce the jump and the average operators, which are
commonly used in the discontinuous Galerkin framework. 
Let $e \in \MEhI$ be any interior face shared by two neighbouring
elements $K^+$ and $K^-$, i.e. $e = \partial K^+ \cap \partial K^-$,
and we let $\un^+$ and $\un^-$ be the unit outward normal vectors on
$e$ corresponding to $K^+$ and $K^-$, respectively. 
Let $v$ be a piecewise smooth scalar-valued function, and let $\bm{q}$
be a piecewise smooth vector- or tensor-valued function, 
the jump operator $\jump{\cdot}$ and the average $\aver{\cdot}$ on $e$
are defined as 
\begin{displaymath}
  \jump{v}|_e := v^+|_e - v^-|_e, \quad \aver{v}|_e =
  \frac{1}{2}(v^+ |_e+ v^-|_e), \quad \jump{\bm{q}}|_e := \bm{q}^+|_e
  - \bm{q}^-|_e, \quad \aver{\bm{q}}|_e =  \frac{1}{2}(\bm{q}^+ |_e+
  \bm{q}^-|_e), 
\end{displaymath}
where $v^{\pm} := v|_{K^{\pm}}$, $\bm{q}^{\pm} := \bm{q}|_{K^{\pm}}$.
On the boundary face $e \in \MEhB$, both operators on $e$ are modified
as 
\begin{displaymath}
  \jump{v}|_e := \aver{v}|_e := \tilde{v}|_e, \quad 
  \jump{\bm{q}}|_e := \aver{\bm{q}}|_e := \tilde{\bm{q}}|_e,
\end{displaymath}
where $\tilde{v} := v|_K$, $\tilde{\bm{q}} := \bm{q}|_K$ with $e
\subset \partial K$.

Given an open bounded domain $D \subset \mb{R}^d$, we let $H^s(D)$ 
denote the usual Sobolev space with the exponent $s \geq 0$, and its
associated inner products, seminorms and norms are also followed.
Throughout this paper, $C$ and $C$ with subscripts are denoted as 
generic constants that may vary in different lines but are always
independent of $h$.

Our model problem is the following biharmonic equation, which reads
\begin{equation}
  \begin{aligned}
    \Delta^2 u & = f, && \text{in } \Omega, \\
    u & = 0, && \text{on } \partial \Omega, \\
    \partial_{\un} u & = 0, && \text{on } \partial \Omega, 
  \end{aligned}
  \label{eq_biharmonic}
\end{equation}
where $f \in L^2(\Omega)$, and $\un$ denotes the unit outward normal
vector on the boundary $\partial \Omega$. 
For simplicity, we assume that the problem \eqref{eq_biharmonic}
possesses a unique solution $u \in H^4(\Omega)$. 
We refer to \cite{Blum1980boundary} for more results about the
regularity. 


\section{Reconstruction from $C^0$ linear space}
\label{sec_space}
In this section, we introduce an approximation space for the problem
\eqref{eq_biharmonic} by presenting a reconstruction procedure 
on the standard $\mP_1$ finite element spaces $V_h := \{ v_h \in
H^1(\Omega): \ v_h|_K \in \mP_1(K), \ \forall K \in \MTh\}$ and $V_{h,
0} := V_h \cap H_0^1(\Omega)$.
The main idea is to define a linear operator $\mR^m$ that maps any
$v_h \in V_h$ into a high-order piecewise polynomial function, i.e.
its image space will be a piecewise polynomial space.
The reconstruction contains two steps: constructing a wide element
patch for each element and solving a local least squares fitting
problem on each patch.

\textbf{Step 1.}  
For each $K \in \MTh$, we collect some neighbouring elements to form a
patch set $\mS(K)$, which is carried out by a recursive algorithm. 
We begin by assigning a constant $N_{m}$ to govern the size of the
patch $\mS(K)$. 
The value of $N_m$ will be specified later.
We define a sequence of sets $\mS_t(K)$ in a recursive manner, which
read
\begin{equation}
  \mS_0(K) := \{K \}, \quad \mS_{t+1}(K) := 
  \bigcup_{K' \in \mS_t(K)} w(K'), \quad t \geq 0.
  \label{eq_Sk}
\end{equation}
For each $\mS_t(K)$, we define a corresponding set $\mI_t(K) :=
\bigcup_{K' \in S_t(K)} \mN_{K'}$, which consists of vertices to all
elements in $\mS_t(K)$. 
The recursive algorithm \eqref{eq_Sk} stops when the depth
$t$ meets the condition $ \# \mc{I}_t(K) < N_{m} \leq \#
\mc{I}_{t+1}(K)$, and meanwhile the recursive depth is denoted as
$t_K$.
Then, we let $\mS(K) := \mS_t(K)$ and sort the elements in 
$\mS_{t+1}(K) \backslash \mS_t(K)$ by their distances to $K$, where 
the distance between two elements is measured by the length of the
line connecting their barycenters.
After sorting, we add the first $n_t$ elements in $\mS_{t+1}(K)
\backslash \mS_t(K)$ to the set $\mS(K)$  
such that $\mI(K) := \bigcup_{K' \in \mS(K)} \mN_{K'}$ has at least
$N_{m}$ nodes.  
The detailed algorithm to construct $\mS(K)$ and $\mI(K)$ is presented
in Algorithm~\ref{alg_patch}.
For the set $\mS(K)$, we define its corresponding domain
$\mD(K) := \text{Int}(\bigcup_{K' \in \mS(K)}
\overline{K'})$.  
By the quasi-uniformity of $\MTh$, there holds $h_{\mD(K)} \leq C t_K
h_K$, where $h_{\mD(K)} := \text{diam}(\mD(K))$.

\begin{algorithm}[t]
  \caption{Construction to $\mS(K)$ and $\mI(K)$;}
  \label{alg_patch}
  \KwIn{The mesh $\MTh$ and the threshold $N_{m}$;}
  \KwOut{$\mS(K)$ and $\mI(K)$ for all elements;}
  \For{$K \in \MTh$ }{set $t = 0$, $\mS_0(K) = \{K\}$, $\mI_0(K) =
  \mN_{K}$; \\
  \Repeat{ $\# \mI_t(K) > N_{\mI}$}{
  let $\mS(K) = \mS_t(K)$; \\
  construct $\mS_{t+1}(K)$ by \eqref{eq_Sk};  \\
  let $\mI_{t+1}(K) = \bigcup_{K' \in \mS_{t+1}(K)} \mN_{K'}$; \\
  $t = t+1$;
  }
  let $\tilde{\mS}(K) = \mS_t(K) \backslash \mS(K)$ and sort the elements
  in $\tilde{\mS}(K)$ by their distances to $K$; \\
  \For{$\tilde{K} \in \tilde{\mS}(K)$}{ 
  add $\tilde{K}$ to $\mS(K)$; \\ 
  let $\mI(K) = \bigcup_{K' \in \mS(K)} \mN_{K'}$; \\
  \If{$\# \mI(K) \geq N_{m}$}{
  return $\mS(K)$ and $\mI(K)$;}}}
\end{algorithm}

\textbf{Step 2.} For each $K \in \MTh$, we solve a local least squares
fitting problem based on nodes in $\mI(K)$ and the space $V_h$.
Given any $v_h \in V_h$, we seek a polynomial of degree $m$ by the
following constrained least squares problem: 
\begin{equation}
  \begin{aligned}
    \mathop{\arg\min}_{p \in \mP_m(\mD(K))}  \sum_{\bm{x} \in
    \mI(K)}& | p(\bm{x}) - v_h(\bm{x})|^2,   \\
     \text{s.t. } p(\bm{y}) = v_h(&\bm{y}), \quad \forall \bm{y} \in
     \mN_{K}.
  \end{aligned}
  \label{eq_leastsquares}
\end{equation}

The unisolvence of the problem \eqref{eq_leastsquares} 
depends on the distribution of nodes in $\mI(K)$, which are required
to deviate from being located on an algebraic curve of degree $m$. 
In the following assumption, we give an equivalent statement of the
unisolvence:
\begin{assumption}
  For any $K \in \MTh$ and any $p \in
  \mP_m(\mD(K))$, $p|_{\mI(K)} = 0$ implies $p|_{\mD(K)} = 0$.
  \label{as_1}
\end{assumption}
A necessary condition to Assumption \ref{as_1} is that $\# \mI(K) \geq
\dim(\mP_m(\mD(K)))$, which can be easily fulfilled by choosing a bit
large $N_{m}$. 
The existence and the uniqueness of the solution to
\eqref{eq_leastsquares} are established in the following lemma. 
\begin{lemma}
  For each $K \in \MTh$, the problem \eqref{eq_leastsquares} admits a
  unique solution.
  \label{le_unique_sol}
\end{lemma}
\begin{proof}
  we mainly prove the uniqueness of the solution to
  \eqref{eq_leastsquares}.
  Let $I_K: C^0(K) \rightarrow \mP_1(K)$ be the standard Lagrange
  interpolation operator on $K$.  
  Let $p$ be any solution in \eqref{eq_leastsquares}, and clearly $p +
  t(q - I_K q)$ satisfies the constraint in \eqref{eq_leastsquares}
  for any $q \in \mP_m(\mD(K))$ and any $t \in \mb{R}$.
  We derive that 
  \begin{displaymath}
    \sum_{\bm{x} \in \mI(K)} | p(\bm{x}) - v_h(\bm{x})|^2 \leq
    \sum_{\bm{x} \in \mI(K)} | p(\bm{x}) + t (q - I_Kq)(\bm{x}) -
    v_h(\bm{x})|^2,
  \end{displaymath}
  which leads to 
  \begin{displaymath}
    2t \sum_{\bm{x} \in \mI(K)} (q - I_K q)(\bm{x}) \cdot (p -
    v_h)(\bm{x}) + t^2 ( (q - I_K q)(\bm{x}))^2 \geq 0, \quad \forall
    q \in \mP_m(\mD(K)), \quad \forall t \in \mb{R}.
  \end{displaymath}
  Since $t$ is arbitrary, there holds 
  \begin{equation}
    \sum_{\bm{x} \in \mI(K)} (q - I_K q)(\bm{x}) \cdot (p -
    v_h)(\bm{x}) = 0, \quad \forall q \in \mP_m(\mD(K)).
    \label{eq_qpiuorth}
  \end{equation}
  Let $p_1, p_2$ be the solutions to \eqref{eq_leastsquares}, and we
  know that $p_1 - p_2 = p_1 - p_2 - I_K(p_1 - p_2)$. 
  Bringing $q = p_1 - p_2$ into \eqref{eq_qpiuorth} yields that $p_1 -
  p_2$ vanishes at all nodes in $\mI(K)$. 
  Then, Assumption \ref{as_1} indicates $p_1 = p_2$, which
  immediately confirms that 
  the problem \eqref{eq_leastsquares} admits a unique solution. 
  This completes the proof.
\end{proof}
It is noticeable that the solution of the least squares problem
\eqref{eq_leastsquares} has a linear dependence on the given function
$v_h$. 
This fact inspires us to define a linear map $\mR_K^m: V_h \rightarrow
\mP_m(\mD(K))$ such that $\mR_K^m v_h$ is the solution of the problem
\eqref{eq_leastsquares} for any $v_h \in V_h$.
From the local operator $\mR_K^m$, it is natural to define a global
operator $\mR^m$ in an elementwise manner, which reads 
\begin{equation}
  \begin{aligned}
    \mR^m : \,\,  V_h &\longrightarrow \mV_h^m, \\
    v_h &\longrightarrow \mR^m v_h,
  \end{aligned}
  \quad (\mR^m v_h)|_K := (\mR^m_K v_h)|_K, \quad \forall K \in
  \MTh.
  \label{eq_Rm}
\end{equation}
Here $\mV_h^m := \mR^m V_h$ is the image space of the operator
$\mR^m$.
By the definition \eqref{eq_Rm}, we have that $\mR^m v_h(\forall v_h
\in V_h)$ is a piecewise polynomial function of degree $m$ and
involves the discontinuity across the interelement faces, i.e.
$\mV_h^m$ is a discontinuous piecewise polynomial space. 

Next, we present some properties of $\mR^m$ and $\mV_h^m$.  
The first conclusion is that the operator $\mR^m$ is full-rank.
\begin{lemma}
  The operator $\mR^m$ is non-degenerate, and $\dim(\mV_h^m) =
  \dim(V_h)$.
  \label{le_fullrank}
\end{lemma}
\begin{proof}
  From the linearity of $\mR^m$, it suffices to prove that if any
  function $v_h \in V_h$ satisfies that $\mR^m v_h = 0$, then $v_h$
  must be the zero function. 
  It is evident that $\mR^m v_h = 0$ gives that $\mR^m_K v_h = 0$ for
  $\forall K \in \MTh$. 
  By the constraint in \eqref{eq_leastsquares}, there holds $(\mR^m_K
  v_h)(\bm{y}) = v_h(\bm{y}) = 0$ for $\forall \bm{y} \in \mN_K$,
  which directly implies that $v_h|_K = 0$ for $\forall K \in \MTh$. 
  Hence, we conclude that $\mR^m$ is full-rank, which also brings us
  the result $\dim(\mV_h^m) = \dim(V_h)$.
  This completes the proof.
\end{proof}
Lemma \ref{le_fullrank} is essentially based on the constraint in
\eqref{eq_leastsquares}. 
We note that the constraint in \eqref{eq_leastsquares} is fundamental
in our method, which provides the non-degenerate property of the
reconstruction operator and further allows us to develop the
preconditioning method.
Given any $w_h \in \mV_h^m$, again by the constraint in
\eqref{eq_leastsquares} the function $v_h$ can be determined by
$v_h|_K = I_K (w_h|_K)$ for $\forall K \in \MTh$ such that $\mR^m v_h
= w_h$.
Hereafter, we extend the interpolant polynomial $I_K v(\forall v \in
C^0(K))$ to the domain $\mD(K)$ by the direct polynomial extension,
i.e. $I_K v \in \mP_1(\mD(K))$.
Since $I_K v$ is linear, we derive that
\begin{equation}
  \begin{aligned}
    \| \nabla (I_K v) \|_{L^2(\mD(K))} & \leq C (h_{\mD(K)} /
    h_K)^{d/2} \| \nabla (I_K v) \|_{L^2(K)} \leq C t_K^{d/2}  \|
    \nabla (I_K v) \|_{L^2(K)}, \\
    \| I_K v \|_{L^{\infty}(\mD(K))} & \leq \| I_K v
    \|_{L^{\infty}(K)} + h_{\mD(K)} \| \nabla (I_K v)
    \|_{L^{\infty}(K)} \leq C t_K \| I_K v \|_{L^{\infty}(K)}, 
  \end{aligned} \quad \forall v \in C^0(K).
  \label{eq_IKv}
\end{equation}
Moreover, we outline a group of basis functions to the space
$\mV_h^m$.
Let $\omega_{\bm{\nu}}$ denote the Lagrange basis funtcion with
respect to the node $\bm{\nu} \in \MNh$, i.e. $V_h =
\text{span}(\{\omega_{\bm{\nu}} \}_{\bm{\nu} \in \MNh})$.
Because $\mR^m$ is invertible, we let $\varphi_{\bm{\nu}} := \mR^m
\omega_{\bm{\nu}}$ for $\forall \bm{\nu} \in \MNh$, and then $ \{
\varphi_{\bm{\nu}} \}$ are basis functions of $\mV_h^m$, i.e.
$\mV_h^m = \text{span}( \{ \varphi_{\bm{\nu}}
\}_{\bm{\nu} \in \MNh})$. 
For any $\bm{\nu} \in \MNh$, there holds
$\omega_{\bm{\nu}}(\bm{\nu}) = 1$ while $\omega_{\bm{\nu}}$ vanishes
at all other nodes. This property indicates 
$\omega_{\bm{\nu}}|_{\mD(K')} = 0$ for all $K'$ satisfying 
$\bm{\nu} \not\in \mI(K')$, which
implies that $\varphi_{\bm{\nu}}$ is compactly supported with 
$\text{supp}(\varphi_{\bm{\nu}}) = \bigcup_{K': \  \bm{\nu} \in
\mI(K')} \overline{K'}$. 
By the group of $\{ \varphi_{\bm{\nu}} \}_{\bm{\nu} \in \MNh}$, any
$\mR^m v_h$ can be expanded as 
\begin{equation}
  \mR^m v_h = \sum_{\bm{\nu} \in \MNh} v_h(\bm{\nu})
  \varphi_{\bm{\nu}}, \quad \forall v_h \in V_h.
  \label{eq_Rmexpand}
\end{equation}
The expansion \eqref{eq_Rmexpand} can be directly extended to the
space $C^0(\Omega)$. 
For any $v \in C^0(\Omega)$, we define $\mR^m v$ as 
\eqref{eq_Rmexpand}, or, equivalently, we define $\mR^m v := \mR^m
v_h$, where $v_h$ is the interpolant of $v$ into the space $V_h$.

Let us focus on the stability property on $\mR^m$. 
We introduce the following constants, which measure the stability of
the operator $\mR^m$ in some sense. 
\begin{equation}
  \Lambda_m := \max_{K \in \MTh} (1 + \Lambda_{m, K} t_K \sqrt{\#
  \mI(K)}), \quad 
  \Lambda_{m, K}^2 := \max_{p \in \mP_m(\mD(K))} \frac{ \| p
  \|_{L^2(K)}^2 }{ h_K^d \sum_{\bm{x} \in \mI(K)}  p(\bm{x})^2}, \quad
  \forall K \in \MTh.
  \label{eq_LambdamK}
\end{equation}
The stability estimate is presented in the following lemma. 
\begin{lemma}
  There holds 
  \begin{equation}
    \| \mR^m_K v \|_{L^2(K)} \leq C \Lambda_m h_K^{d/2}
    \max_{\bm{x} \in I(K)} |v|, \quad \forall K \in \MTh, \quad
    \forall v \in C^0(\Omega).
    \label{eq_recon_stability}
  \end{equation}
  \label{le_recon_stability}
\end{lemma}
\begin{proof}
  Let $p = \mR_K^m v$ be the solution to \eqref{eq_leastsquares}.
  By setting $q = p$ in \eqref{eq_qpiuorth}, 
  it can be seen that 
  \begin{displaymath}
    \sum_{\bm{x} \in \mI(K)} (p - I_K p)(\bm{x}) \cdot (p -
    v)(\bm{x}) = 0.
  \end{displaymath}
  This orthogonal property indicates that 
  \begin{displaymath}
    \sum_{\bm{x} \in \mI(K)} (( p - I_K p)(\bm{x}))^2 \leq
    \sum_{\bm{x} \in \mI(K)} ((v - I_K p)(\bm{x}))^2.
  \end{displaymath}
  Since $I_K p = I_K v$, and by the definition \eqref{eq_LambdamK} and
  the estimate \eqref{eq_IKv}, we
  know that
  \begin{align*}
    \| p - I_K v \|_{L^2(K)}^2 & \leq \Lambda_{m, K}^2 h_K^d
    \sum_{\bm{x} \in \mI(K)}  (( p - I_K v )(\bm{x}))^2 \leq
    \Lambda_{m, K}^2 h_K^d \sum_{\bm{x} \in \mI(K)}  (( v - I_K
    v)(\bm{x}))^2 \\
    & \leq C \Lambda_{m, K}^2 h_K^d  \# \mI(K) ( 1 + t_K^2)
    \max_{\bm{x} \in \mI(K)} |v(\bm{x})|^2,
  \end{align*}
  and 
  \begin{displaymath}
    \|I_K v\|_{L^2(K)}^2 \leq C h_K^d \max_{\bm{\nu} \in \mN_{K}}
    |v(\bm{\nu})|^2 \leq  C h_K^d \max_{\bm{x} \in \mI(K)}
    |v(\bm{x})|^2. 
  \end{displaymath}
  Combining the above estimates leads to \eqref{eq_recon_stability},
  which completes the proof.
\end{proof}
By \eqref{eq_recon_stability}, it is quite formal to give the
following approximation error estimate, and we refer to \cite[Theorem
3.3]{Li2012efficient} for the proof.
\begin{lemma}
  There exists a constant $C$ such that
  \begin{equation}
    \| v - \mR_K^m v \|_{H^q(K)} \leq C \Lambda_m h_K^{s - q} \| v
    \|_{H^{s}(\mD(K))}, \quad \forall K \in \MTh, \quad \forall v \in
    H^s(\Omega).
    \label{eq_Rmapp}
  \end{equation}
  where $0 \leq q \leq s - 1$ and $2 \leq s \leq m + 1$.
  \label{le_Rmapp}
\end{lemma}
From Lemma \ref{le_Rmapp}, we find that $\mR_K^m$ has the optimal
convergence rate if $\Lambda_m$ admits a upper bound independent of
$h$. 
Generally speaking, this condition can be fulfilled by selecting a
large enough element patch. 
More details about $\Lambda_m$ and $\Lambda_{m, K}$ are presented in
Remark \ref{re_lam}.
It is noted that the constant $\Lambda_{m, K}$ indeed corresponds to
the minimum singular value of a local matrix, which can be easily
computed, see Appendix \ref{sec_reconM}.
In the appendix, we also present a series of numerical tests on
$\Lambda_m$. 
It can be observed that $\Lambda_m$ has a uniform upper
bound by selecting a large threshold $N_m$.
\begin{remark}
  The constant $\Lambda_{m, K}$ is related to the constant
  \begin{displaymath}
    \Theta_{m, K} := \max_{p \in \mP_m(\mD(K))} \frac{\max_{\bm{x}
    \in \mD(K)} p(\bm{x}) }{\max_{\bm{x} \in \mI(K)} p(\bm{x}) },
    \quad \forall K \in \MTh,
  \end{displaymath}
  and it is clear that $C \Theta_{m, K} \leq \Lambda_{m, K} \leq
  \Theta_{m, K}$. 
  Here $\Theta_{m, K}$ is close to the Lebesgue constant
  \cite{Powell1981approximation}, and, unfortunately, to our best
  knowledge, there are few results to the upper bound of the Lebesgue
  constant in two and three dimensions.
  Currently, we can prove that for a wide enough element patch
  $\mD(K)$, $\Theta_{m, K}$ admits a uniform upper bound. 
  Let $B_{r, K}$ and $B_{R, K}$ be the largest and the smallest balls
  such that $B_{r, K} \subset \mD(K) \subset B_{R, K}$ with the
  radii $r_K$ and $R_K$, respectively. 
  By \cite[Lemma 5]{Li2016discontinuous}, we have that $\Theta_{m, K}
  = 2$ under the condition that $r_K \geq 2m \sqrt{R_K h_K}$. 
  Generally speaking, there will be $R_K \approx r_K \approx O(t_K
  h_K)$ for the wide enough patch. Thus, this condition can be
  fulfilled when $R_K \approx r_K$.
  In \cite[Lemma 6]{Li2016discontinuous} and 
  \cite[Lemma 3.4]{Li2012efficient}, we show that there exists a
  threshold $N_{\mS}(m, C_\nu)$ that  only depends on $m, C_\nu$ such
  that $r_K \geq 2m \sqrt{R_K h_K}$ is met when $\# \mS(K) \geq
  N_{\mS}$. 
  Consequently, $\Lambda_m \leq \Theta_m$ also admits a uniform upper
  bound under this condition.
  But $N_{\mS}(m, C_\nu)$ is usually too large and
  impractical in the computer implementation. 
  In this paper, the constant $\Lambda_{m, K}$ can be directly
  computed and the value can serve as an indicator to show whether
  the set $\mI(K)$ is proper.
  In Appendix \ref{sec_reconM}, we present the method to compute
  $\Lambda_m$ for a given mesh. 
  From the numerical tests, we observe that the condition that
  $\Lambda_m$ admits a upper bound can be met by choosing a large
  $N_m$. 
  Roughly speaking, the threshold $N_m$ is approximately equal to
  $1.5\dim(\mP_m(\cdot))$.
  \label{re_lam}
\end{remark}

To end this section, we define the approximation space $U_h^m$ that
will be used to approximate the solution of the problem
\eqref{eq_biharmonic} in next section. 
Here $U_h^m$ is defined as $U_h^m := \mR^m V_{h, 0}$. 
Since $V_{h, 0} = V_h \cap H_0^1(\Omega)$, and by the properties of
$\mR^m$, it is similar to verify that $U_h^m$ has the following
properties:
\begin{enumerate}
  \item[1.] $\dim(U_h^m) = \dim(V_{h, 0})$ and $U_h^m = \text{span}(
    \{ \varphi_{\bm{\nu}} \}_{\bm{\nu} \in \MNhI})$.
  \item[2.] For any $v \in H^s(\Omega) \cap H_0^1(\Omega)$, there
    holds $\mR^m v \in U_h^m$ satisfying the approximation estimate
    \eqref{eq_Rmapp} in each element.
\end{enumerate}
Compared with the space $\mV_h^m$, $U_h^m$ has less degrees of freedom
and it will also provide a simpler preconditioning method for us.


\section{Numerical Scheme} 
\label{sec_scheme}
In this section, we present and analyze the numerical scheme for
solving the biharmonic equation \eqref{eq_biharmonic}, based on the
space $U_h^m$.  
Since $U_h^m$ is a discontinuous piecewise polynomial space, we are
allowed to adopt the symmetric interior penalty
discontinuous Galerkin scheme \cite{Mozolevski2007hp,
Georgoulis2008discontinuous} to seek the numerical solution. 
The discrete variational problem reads: seek $u_h \in U_h^m$ such
that 
\begin{equation}
  a_h(u_h, v_h) = l_h(v_h), \quad \forall v_h \in U_h^m,
  \label{eq_dvar}
\end{equation}
where 
\begin{equation}
  \begin{aligned}
    a_h(v_h, w_h) &:= \sum_{K \in \MTh} \int_K \lap v_h \lap w_h
    \dx{x} + \sum_{e \in \MEh} \int_e \left(\jump{v_h} 
    \aver{\nabla_{\un} \lap w_h} + \jump{w_h} \aver{\nabla_{\un}
    \lap v_h} \right) \dx{s} \\
    &- \sum_{e \in \MEh} \int_e \left( \aver{\lap w_h} \jump{
    \nabla_{\un} v_h} + \aver{\lap v_h} \jump{\nabla_{\un} w_h} \right)
    \dx{s} \\
    &+ \sum_{e \in \MEh} \int_e \left( \mu_1 h_e^{-3} \jump{v_h} 
    \jump{w_h} \dx{s} + \mu_2 h_e^{-1} \jump{\nabla_{\un} v_h} \jump{
    \nabla_{\un} w_h} \right) \dx{s}, \quad \forall v_h, w_h \in U_h,
  \end{aligned}
  \label{eq_ah}
\end{equation}
and 
\begin{displaymath}
  l_h(v_h) := \sum_{K \in \MTh} \int_K f v_h \dx{x}, \quad \forall v_h
  \in U_h.
\end{displaymath}
Here $U_h := U_h^m + H^4(\Omega)$ and $\mu_1$, $\mu_2$ are referred as
penalty parameters.
We refer to \cite{Georgoulis2008discontinuous} for the detailed
derivation of the interior penalty form. 

Because $U_h^m$ is a subspace of the standard
discontinuous piecewise polynomial space, the error
estimation to the problem \eqref{eq_dvar} can be directly derived in
the standard DG framework. 

We introduce the following energy norms: 
\begin{displaymath}
  \begin{aligned}
    \DGenorm{v_h}^2 & := \sum_{K \in \MTh} \| \lap v_h \|^2_{L^2(K)} + 
    \sum_{e \in \MEh} h_e^{-3} \| \jump{v_h} \|^2_{L^2(e)} + \sum_{e
    \in \MEh} h_e^{-1} \| \jump{\nabla_{\un} v_h} \|^2_{L^2(e)}, 
    \quad \forall v_h \in U_h, \\
    \DGtenorm{v_h}^2 & := \DGenorm{v_h}^2 + \sum_{e \in \MEh} h_e^3 \|
    \aver{\nabla_{\un} \lap v_h} \|^2_{L^2(e)} + \sum_{e \in \MEh} h_e
    \| \aver{\lap v_h} \|^2_{L^2(e)}, \quad \forall v_h \in U_h.
  \end{aligned}
\end{displaymath}
Both energy norms are equivalent over the piecewise polynomial spaces,
i.e. 
\begin{equation}
  \DGenorm{v_h} \leq \DGtenorm{v_h} \leq C \DGenorm{v_h}, \quad
  \forall v_h \in U_h^m.
  \label{eq_DGnorm_equi}
\end{equation}
The estimate \eqref{eq_DGnorm_equi} follows from the inverse
estimate and the trace estimate; see \cite{Mozolevski2007hp} for the
proof. 
Furthermore, we give the relationship between the energy norm
$\DGenorm{\cdot}$ and $\| \cdot \|_{L^2(\Omega)}$, which reads
\begin{equation}
  \| v_h \|_{L^2(\Omega)} \leq C \DGenorm{v_h} \leq Ch^{-2}  \| v_h
  \|_{L^2(\Omega)}, \quad \forall v_h \in U_h^m.
  \label{eq_DGL2_relation}
\end{equation}
The second estimate is straightforward from the inverse estimate.
The first estimate in \eqref{eq_DGL2_relation}, i.e. $\DGenorm{\cdot}$
is stronger than $\| \cdot \|_{L^2(\Omega)}$, can be verified by the
dual argument. 
Given any $v_h \in U_h^m$, we consider the elliptic problem 
\begin{displaymath}
  -\Delta w = v_h, \ \text{in } \Omega, \quad w = 0, \ 
  \text{on } \partial \Omega.
\end{displaymath}
Since $\Omega$ is convex, the above problem admits a unique
solution $w \in H^2(\Omega) \cap H_0^1(\Omega)$ with $\|w
\|_{H^2(\Omega)} \leq C \| v_h \|_{L^2(\Omega)}$.
Applying integration by parts and the trace estimate, we find that 
\begin{align*}
  \| v_h \|_{L^2(\Omega)}^2  & = -\int_{\Omega} \Delta w v_h \dx{x} =
  - \sum_{e \in \MEh} \int_e (\nabla_{\un} w) \jump{v_h} \dx{s} +
  \sum_{e \in \MEh} \int_e w \jump{\nabla_{\un} v_h} \dx{s} - \sum_{K
  \in \MTh} \int_K w \Delta v_h \dx{x} \\
  & \leq C ( \| w \|_{L^2(\Omega)}^2 + \sum_{e \in \MEh} (h_e^3 \|
  \nabla_{\un} w \|_{L^2(e)}^2 + h_e \| w \|_{L^2(e)}^2))^{1/2}
  \DGenorm{v_h}  \\
  & \leq C \| w \|_{H^2(\Omega)} \DGenorm{v_h} \leq C \|
  v_h \|_{L^2(\Omega)} \DGenorm{v_h}.
\end{align*}
Thus, the desired estimate \eqref{eq_DGL2_relation} is reached.

In the analysis to the preconditioned system, we need another 
energy norm $\snorm{\cdot}$, which reads
\begin{displaymath}
  \snorm{v_h}^2 := \sum_{K \in \MTh} \| D^2 v_h \|^2_{L^2(K)} + 
  \sum_{e \in \MEh} h_e^{-3} \| \jump{u_h} \|^2_{L^2(e)} + \sum_{e
  \in \MEh} h_e^{-1} \| \jump{\nabla_{\un} v_h} \|^2_{L^2(e)},
  \quad \forall v_h \in U_h.
\end{displaymath}
We note that both energy norms $\DGenorm{\cdot}$ and $\snorm{\cdot}$
are also equivalent restricted on $U_h^m$: 
\begin{equation}
  \DGenorm{v_h} \leq C \snorm{v_h} \leq C \DGenorm{v_h}, \quad \forall
  v_h \in U_h^m.
  \label{eq_DGsnormequ}
\end{equation}

\begin{remark}
  In \eqref{eq_DGsnormequ}, the norm $\snorm{\cdot}$ is stronger than 
  $\DGenorm{\cdot}$ is trivial while the reverse estimate 
  can be obtained by the discrete Miranda-Talenti inequality,
  which reads 
  \begin{equation}
    \sum_{K \in \MTh} \| D^2 v_h \|_{L^2(K)}^2 \leq \sum_{K \in \MTh}
    \| \Delta v_h \|_{L^2(K)}^2 + C \sum_{e \in \MEhI}
    h_e^{-3} \| \jump{v_h} \|_{L^2(e)}^2 + C \sum_{e \in \MEh}
    h_e^{-1} \| \jump{\nabla_{\un} v_h } \|_{L^2(e)}^2,
    \label{eq_dMT}
  \end{equation}
  for $\forall v_h \in U_h^m$.  In \cite{Neilan2019discrete}, the
  authors proved that the estimate \eqref{eq_dMT} holds for 
  the piecewise polynomial space of degree $m$, where $m \geq 2$ in
  two dimensions and $2 \leq m \leq 3$ in three dimensions. 
  In Appendix \ref{sec_app_dMT}, we present another proof to
  \eqref{eq_dMT} for any $m \geq 2$ in both two and three dimensions. 
  From \eqref{eq_dMT}, the equivalence \eqref{eq_DGsnormequ} can be
  easily verified.
  \label{re_DGsnormequ}
\end{remark}

The bilinear form $a_h(\cdot, \cdot)$ is bounded and coercive under
the energy norm $\DGtenorm{\cdot}$, and the Galerkin orthogonality
holds. 
These results are standard in the DG framework, and the detailed
proofs are referred to \cite{Mozolevski2007hp,
Georgoulis2008discontinuous}.

\begin{lemma}
  Let $a_h(\cdot, \cdot)$ be defined with the sufficiently large
  $\mu_1$ and $\mu_2$, there hold
  \begin{align}
    |a_h(v_h, w_h)| & \leq C \DGtenorm{v_h} \DGtenorm{w_h},
    \quad \forall v_h, w_h \in U_h, \label{eq_ahb} \\
    a_{h}(v_h, v_h) & \geq C \DGtenorm{v_h}^2, \quad \forall
    v_h \in U_h^m.
    \label{eq_ahc}
  \end{align}
  \label{le_ahbc}
\end{lemma}
\begin{lemma}
  Let $u \in H^4(\Omega)$ be the exact solution to the problem
  \eqref{eq_biharmonic}, and let $u_h \in U_h^m$ be the numerical
  solution to \eqref{eq_dvar}, then there holds
  \begin{equation}
    a_h(u - u_h, v_h) = 0, \quad \forall v_h \in U_h^m.
    \label{eq_ahorth}
  \end{equation}
  \label{le_ahorth}
\end{lemma}
The convergence analysis follows from Lemma \ref{le_ahbc} - Lemma
\ref{le_ahorth} and the approximation property of the space $U_h^m$.
\begin{theorem}
  Let $u \in H^{t}(\Omega)(t \geq 4)$ be the exact solution to
  \eqref{eq_biharmonic}, and let $u_h \in U_h^m$ be the numerical
  solution to \eqref{eq_dvar} with $m \geq 2$, and
  let $\mu_1, \mu_2$ be taken as in Lemma \ref{le_ahbc}, then there
  holds
  \begin{equation}
    \DGenorm{u - u_h} \leq C \Lambda_m h^{s - 2} \| u
    \|_{H^{t}(\Omega)}, \quad s = \min(m + 1, t).
    \label{eq_DGerror}
  \end{equation}
  \label{th_DGerror}
\end{theorem}
\begin{proof}
  Let $u_I := \mR^m u$ be the interpolant of $u$ into the space
  $U_h^m$. 
  Combining the approximation estimate \eqref{eq_Rmapp} and
  the trace estimate, we obtain that
  \begin{displaymath}
    \DGtenorm{u - u_I} \leq C \Lambda_m h^{s-2} \| u \|_{H^t(\Omega)}.
  \end{displaymath}
  By \eqref{eq_ahb} - \eqref{eq_ahorth}, we have that
  \begin{align*}
    \DGtenorm{u_I - u_h}^2 & \leq C a_h(u_I - u_h, u_I - u_h) = C
    a_h(u_I - u, u_I - u_h) \leq C \DGtenorm{u_I - u} \DGtenorm{u_I -
    u_h}, 
  \end{align*}
  and by the triangle inequality, we obtain that
  \begin{displaymath}
    \DGtenorm{u - u_h} \leq C \DGtenorm{u - u_I} \leq C \Lambda_m
    h^{s-2} \| u \|_{H^t(\Omega)},
  \end{displaymath}
  which completes the proof.
\end{proof}
The $L^2$ error estimate is obtained by the duality argument.  
Assume that the dual problem 
\begin{equation}
  \begin{aligned}
    \lap^2 \psi &=  u - u_h,  && \text{in } \Omega, \\
    \psi &= 0, && \text{on } \partial \Omega, \\
    \partial_{\un} \psi &= 0, && \text{on } \partial \Omega,
  \end{aligned}
  \label{eq_duality}
\end{equation}
admits a unique solution $\psi \in H^4(\Omega)$ such that $\|
\psi \|_{H^4(\Omega)} \leq C \| u - u_h \|_{L^2(\Omega)}$. 
Let $\psi_I := \mR^m \psi$, we have that 
\begin{equation}
  \DGtenorm{ \psi - \psi} \leq C \Lambda_m h^s \| \psi
  \|_{H^4(\Omega)},
  \label{eq_dualapp}
\end{equation}
where $s = \min(2, m - 1)$.  Multiplying $u - u_h$ in
\eqref{eq_duality} yields that
\begin{align*}
  \| u - u_h \|_{L^2(\Omega)}^2 & =  \int_{\Omega} \Delta^2 \psi (u -
  u_h) \dx{x} = a_h(\psi, u - u_h) = a_h(\psi - \psi_I, u -
  u_h) \\ 
  & \leq C \DGtenorm{\psi - \psi_I} \DGtenorm{u - u_h} \leq C
  h^{m+s - 1} \| \psi \|_{H^4(\Omega)} \| u \|_{H^t(\Omega)}, 
\end{align*}
By eliminating $  \| u - u_h \|_{L^2(\Omega)}$, we summarize the $L^2$
error estimate as 
\begin{equation}
  \| u - u_h \|_{L^2(\Omega)} \leq C \Lambda_m h^s \| u
  \|_{H^t(\Omega)},  \quad s = \begin{cases}
    2, & m = 2, \\
    \min(m + 1, t), & m \geq 3. \\
  \end{cases}
  \label{eq_L2estimate}
\end{equation}

For the inhomogeneous boundary conditions that $u = g_1$ and
$\nabla_{\un} u = g_2$ on $\partial \Omega$ in the problem
\eqref{eq_biharmonic}, the proposed numerical scheme can be simply
modified to this case. 
Let $g_{1, I} \in V_h$ be defined by 
$g_{1, I}(\bm{\nu}) = g_1(\bm{\nu})$ for $\forall \bm{\nu} \in \MNhB$
and $g_{1, I}(\bm{\nu}) = 0$ for $\forall \bm{\nu} \in \MNhI$. 
We let $w_{1, I} := \mR^m g_{1, I}$.
The discrete variational problem then reads: seek $w_h \in U_h^m$ such
that 
\begin{displaymath}
  a_h(w_h, v_h) = \tilde{l}_h(v_h), \quad \forall v_h \in U_h^m,
\end{displaymath}
where $a_h(\cdot, \cdot)$ is the same as \eqref{eq_ah}, and
$\tilde{l}_h(\cdot)$ is defined as 
\begin{align*}
  \tilde{l}_h(v_h) =  l_h(v_h) - a_h(w_{1, I}, v_h) 
  -& \sum_{e \in \MEhB} \int_e g_2 \aver{\Delta v_h} \dx{s} \\
  + & \sum_{e \in \MEhB} \int_e ( \mu_1 h_e^{-3} g_1 \jump{v_h} + \mu_2
  h_e^{-1} g_2 \jump{\nabla_{\un} v_h} ) \dx{s}, \quad \forall v_h \in
  U_h^m.
\end{align*}
It is similar to show that the solution $u_h := w_h + w_{1, I}$
satisfies the error estimates \eqref{eq_DGerror} and
\eqref{eq_L2estimate}.
Now, we have proved that
our scheme also has the optimal convergence rates as in the standard
discontinuous Galerkin method
\cite{Georgoulis2008discontinuous, Mozolevski2007hp}.

In the rest of this section, we focus on the resulting linear system. 
Generally, the condition number for the linear system grows very fast
at the speed  $O(h^{-4})$ in solving the fourth-order problem.
We first demonstrate that in our scheme the condition number still has
a similar estimate.

Let $A_m$ be the matrix with respect to the bilinear form
$a_h(\cdot, \cdot)$ over spaces $U_h^m \times U_h^m$. 
As stated earlier, 
$U_h^m = \text{span}\{ \varphi_{\bm{\nu}} \}_{\bm{\nu} \in \MNhI}$,
and thus $A_m$ can be expressed as $A_m = (a_h(\varphi_{\bm{\nu}_1},
\varphi_{\bm{\nu}_2}))_{n_p \times n_p}$, where $n_p := \# \MNhI$.
We define the mass matrix  $M_m := ((\varphi_{\bm{\nu}_1},
\varphi_{\bm{\nu}_2})_{L^2(\Omega)})_{n_p \times n_p}$.  
For any finite element function $v_h \in U_h^m$, we associate it with
a coefficient vector $\bmr{v} = \{ v_{\bm{\nu}} \}_{\bm{\nu} \in
\MNhI} \in \mb{R}^{n_p}$, where 
$v_{\bm{\nu}} = v_h(\bm{\nu}) $ for any $\bm{\nu} \in \MNhI$. 
Conversely, any vector $\bmr{v} \in \mb{R}^{n_p}$ also corresponds to
a finite element function $v_h \in U_h^m$ by $v_h(\bm{\nu}) =
v_{\bm{\nu}}(\forall \bm{\nu} \in \MNhI)$.
The linear system $a_h(\cdot, \cdot)$ and the $L^2$ inner product
$(\cdot, \cdot)_{L^2(\Omega)}$ can
be rewritten as 
\begin{displaymath}
  a_h(v_h, w_h) = \bmr{v}^T A_m \bmr{w}, \quad (v_h,
  w_h)_{L^2(\Omega)} = \bmr{v}^T M_m  \bmr{w}, \quad \forall \bmr{v},
  \bmr{w} \in \mb{R}^{n_p},
\end{displaymath}
which directly brings us that
\begin{equation}
  \frac{\bmr{v}^T A_m \bmr{v}}{\bmr{v}^T \bmr{v}} =
  \frac{a_h(v_h,v_h)}{(v_h,v_h)_{L^2(\Omega)}} \frac{\bmr{v}^T M_m
  \bmr{v}}{\bmr{v}^T \bmr{v}}, \quad \forall \bmr{v}\in
  \mb{R}^{n_p} \backslash \{\bm{0}\}.
  \label{eq_kappaA0}
\end{equation}
By Lemma \ref{le_ahbc} and the estimate \eqref{eq_DGL2_relation}, we
conclude that
\begin{displaymath}
  \| v_h \|_{L^2(\Omega)}^2 \leq C \DGenorm{v_h}^2 \leq C a_h(v_h,
  v_h) \leq C \DGenorm{v_h}^2 \leq  Ch^{-4}  \| v_h
  \|_{L^2(\Omega)}^2, \quad \forall v_h \in U_h^m.
\end{displaymath}
It remains to bound the term $(\bmr{v}^T M_m \bmr{v})/(\bmr{v}^T
\bmr{v})$. 
By the inverse estimate, we have that 
\begin{equation}
  \begin{aligned}
    \bmr{v}^T \bmr{v}  & = \sum_{\bm{\nu} \in \MNhI} v_{\bm{\nu}}^2
    \leq C \sum_{K \in \MTh} \sum_{\bm{\nu} \in \mN_K}
    v_{\bm{\nu}}^2 \leq   C \sum_{K \in \MTh} \sum_{\bm{\nu} \in
    \mN_K} (v_h(\bm{\nu}))^2 \\
    & \leq C \sum_{K \in \MTh} \| v_h \|_{L^{\infty}(K)}^2 \leq C
    \sum_{K \in \MTh}
    h^{-d} \| v_h \|_{L^2(K)}^2 = C h^{-d} \bmr{v}^T M_m \bmr{v},
    \quad \forall \bmr{v} \in \mb{R}^{n_p}.
  \end{aligned}
  \label{eq_vTvupper}
\end{equation}
From \eqref{eq_LambdamK}, we have that 
\begin{equation}
  \begin{aligned}
    \bmr{v}^T M_m \bmr{v} = \sum_{K \in \MTh} \| v_h \|_{L^2(K)}^2
    \leq \sum_{K \in \MTh} \Lambda_{m, K}^2 h_K^d \sum_{\bm{\nu} \in
    \mI(K)}  v_{\bm{\nu}}^2 \leq C \Lambda_m^2 h^d \bmr{v}^T \bmr{v},
    \quad \forall \bmr{v} \in \mb{R}^{n_p}.
  \end{aligned}
  \label{eq_vTvlower}
\end{equation}
Combining \eqref{eq_kappaA0} - \eqref{eq_vTvlower} yields that
\begin{displaymath}
  Ch^d \leq \frac{\bmr{v}^T A_m \bmr{v} }{ \bmr{v}^T \bmr{v}} \leq
  C{\Lambda}_m^2 h^{d - 4}, \quad \forall \bmr{v} \in
  \mb{R}^{n_p} \backslash \{ \bm{0} \},
\end{displaymath}
which indicates $\kappa(A_m) \leq C \Lambda_m^2 h^{-4}$.
Therefore, the linear system arising in our method still has the
condition number $O(h^{-4})$, which is highly ill-conditioned as the
mesh is refined.

Now, we present a preconditioning method based on the linear space
$V_{h, 0}$.
It is noticeable that the size of the matrix $A_m$ is always 
$n_p \times n_p$ independent of the degree $m$.
This fact inspires us to construct a preconditioner from the
lowest-order scheme on the same mesh. 
To this end, we introduce the bilinear form $a_h^{\mc{L}}(\cdot,
\cdot)$ by 
\begin{displaymath}
  a_h^{\mc{L}}(v_h, w_h) := \sum_{e \in \MEh}\int_e h_e^{-1} 
  \jump{\nabla_{\un} v_h} \jump{\nabla_{\un} w_h} \dx{s}, \quad
  \forall v_h, w_h \in V_{h, 0}, 
\end{displaymath}
which corresponds to the lowest-order scheme in the sense that 
$a_h^{\mc{L}}(v_h, w_h) = \mu_2 a_h(v_h, w_h)$ for $\forall v_h, w_h
\in V_{h, 0}$. 
An immediate observation is that $ a_h^{\mc{L}}(v_h, v_h) =
\DGenorm{v_h}^2$ for $\forall v_h \in V_{h, 0}$, which implies that 
$a_h^{\mc{L}}(\cdot, \cdot)$ is bounded and coercive
under the energy norm $\DGenorm{\cdot}$. 
Let $A_{\mc{L}}$ be the matrix of $a_h^{\mc{L}}(\cdot, \cdot)$, and
$A_{\mc{L}} \in \mb{R}^{n_p \times n_p}$ is a symmetric positive
definite matrix. 
We will show that $A_{\mc{L}}^{-1}$ can serve as an efficient
preconditioner to the matrix $A_m$ for any accuracy $m$, meaning that
the condition number of the preconditioned system $A_{\mc{L}}^{-1}
A_m$ admits a uniform upper bound independent of $h$.

The main step to estimate the condition number is to establish the
following norm equivalence results for the reconstruction process.
For any $v_h \in V_{h, 0}$, we claim that  $\DGenorm{v_h}$ is
equivalent to $\DGenorm{\mR^m v_h}$ in Lemma \ref{le_DGenorm_upper}
and Lemma \ref{le_DGenorm_lower}.

\begin{lemma}
  There exists a constant $C$ such that 
  \begin{equation}
    \DGenorm{v_h} \leq C \DGenorm{ \mR^m v_h}, \quad \forall v_h
    \in V_{h, 0}.
    \label{eq_DGenorm_upper}
  \end{equation}
  \label{le_DGenorm_upper}
\end{lemma}
\begin{proof}
  Because $v_h$ is piecewise linear, there holds $\DGenorm{v_h}^2 =
  \sum_{e \in \MEh} h_e^{-1} \| \jump{\nabla_{\un} v_h}
  \|_{L^2(e)}^2$. 
  By the equivalence \eqref{eq_DGsnormequ}, it suffices to prove that
  $\DGenorm{v_h} \leq C \snorm{\mR^m v_h}$ for the estimate
  \eqref{eq_DGenorm_upper}.

  For any interior face $e \in \MEhI$ shared by elements $K_1$ and
  $K_2$, we let $\bm{x}_e$ be the barycenter of $e$. 
  Let $v_{K_1} := v_h|_{K_1}, v_{K_2} := v_h|_{K_2}$,
  and we have that  $v_{K_i} = I_{K_i} (\mR_{K_i}^m v_h)(i = 1,2)$.
  Note that $v_{K_i}$ is a linear polynomial, and 
  we apply the inverse estimate and the approximation property of
  $I_{K_i}$ to derive that
  \begin{align*}
    h_e^{-1} \| \jump{\nabla_{\un} v_h} \|^2_{L^2(e)} &\leq C h_e^{d-2} |
    \nabla_{\un}( v_{K_1} - v_{K_2}) |^2 \\
    &\leq C h_e^{d-2} ( | \nabla_{\un}( v_{K_1} -\mR^m_{K_1}
    v_h(\bm{x}_e)) |^2 + |\nabla_{\un} ( (
    \mR^m_{K_2} v_h)(\bm{x}_e) 
    -  v_{K_2}) |^2 ) \\
    & \hspace{130pt} + | \nabla_{\un} ( (\mR^m_{K_1} v_{h})(\bm{x}_e) 
    - ( \mR^m_{K_2} v_{h})(\bm{x}_e)) |^2 \\
    & \leq C h_e^{d-2} ( \| \nabla(\mR_{K_1}^m v_h - I_{K_1}
    \mR_{K_1}^m v_h ) \|_{L^\infty(K_1)}^2 + 
    \| \nabla(\mR_{K_2}^m v_h - I_{K_2}
    \mR_{K_2}^m v_h ) \|_{L^\infty(K_2)}^2) \\
    & \hspace{130pt} + \| \nabla_{\un} (
    \mR_{K_1}^m v_h - \mR_{K_2}^m v_h) \|_{L^{\infty}(e)}^2  \\
    & \leq C ( \| D^2 (\mR_{K_1}^m v_h) \|_{L^2(K_1)}^2 +  \| D^2
    (\mR_{K_2}^m v_h) \|_{L^2(K_2)}^2 + h_e^{-1} \|
    \jump{\nabla_{\un} \mR^m v_h} \|_{L^2(e)}^2).
  \end{align*}
  Summation over all interior faces gives that  $ \sum_{e \in \MEhI}
  h_e^{-1} \| \jump{ \nabla_{\un} v_h } \|_{L^2(e)}^2 \leq C
  \snorm{\mR^m v_h}^2$. 
  For any boundary face $e \in \MEhB$, the same procedure can still be
  applied. 
  For such a face $e$, we let $K$ be the element with $e \subset
  \partial K$.
  We find that 
  \begin{align*}
    h_e^{-1} \| \nabla_{\un} v_h \|_{L^2(e)}^2 & \leq C h_e^{d-2} ( \|
    \nabla_{\un} ( v_h - \mR_K^m v_h) \|_{L^{\infty}(e)}^2 +  \|
    \nabla_{\un} \mR_K^m v_h \|_{L^{\infty}(e)})  \\ 
    & \leq C ( \| D^2 (\mR_K^m v_h) \|_{L^2(K)}^2 + h_e^{-1} \|
    \nabla_{\un} \mR_K^m v_h \|_{L^2(e)}^2),
  \end{align*}
  which directly yields $\sum_{e \in \MEhB} h_e^{-1} \| \jump{
  \nabla_{\un} v_h} \|_{L^2(e)}^2 \leq C \snorm{ \mR^m v_h}^2$.
  Combining above estimates brings us the desired estimate
  \eqref{eq_DGenorm_upper}, which completes the proof.
\end{proof}

\begin{lemma}
  There exists a constant $C$ such that 
  \begin{equation}
    \DGenorm{\mR^m v_h} \leq C \Lambda_m \DGenorm{v_h},
    \quad \forall v_h \in V_{h, 0}.
    \label{eq_DGenorm_lower}
  \end{equation}
  \label{le_DGenorm_lower}
\end{lemma}
\begin{proof}
  The estimation of \eqref{eq_DGenorm_lower} also starts from the
  terms on interior faces. 
  For any $e \in \MEhI$ shared by elements $K_1$ and $K_2$, 
  we let $v_{K_i} := v_h|_{K_i}$, and there holds $v_{K_i} = I_{K_i}
  (\mR_{K_i}^m v_h)$.
  It can be seen that
  \begin{align*}
    h_e^{-1} &\| \jump{\nabla_{\un} \mR^m v_h} \|^2_{L^2(e)} \\
    & \leq
    Ch_e^{-1} ( \| \nabla_{\un} (\mR_{K_1}^m v_h - v_{K_1})
    \|_{L^2(e)} + \| \nabla_{\un} (\mR_{K_2}^m v_h - v_{K_2})
    \|_{L^2(e)} + \| \jump{\nabla_{\un} v_h} \|_{L^2(e)}^2) \\ 
    & \leq C(h_{K_1}^{d-2} \| \nabla (\mR_{K_1}^m v_h - v_{K_1})
    \|_{L^{\infty}(K_1)}^2 + h_{K_2}^{d-2} \| \nabla (\mR_{K_2}^m v_h - v_{K_2})
    \|_{L^{\infty}(K_2 )}^2 + h_e^{-1} \| \jump{\nabla_{\un} v_h}
    \|_{L^2(e)}^2) \\
    & \leq C ( \| D^2 (\mR_{K_1}^m v_h) \|_{L^2(K_1)}^2 +  \| D^2
    (\mR_{K_2}^m v_h) \|_{L^2(K_2)}^2 +  h_e^{-1} \|
    \jump{\nabla_{\un} v_h} \|_{L^2(e)}^2).
  \end{align*}
  Since $v_{K_1}|_e = v_{K_2}|_e$, we derive that 
  \begin{align*}
    h_e^{-3}  \| \jump{\mR^m v_h} \|^2_{L^2(e)}& \leq C h_e^{d - 4} \|
    \mR_{K_1}^m v_h - \mR_{K_2}^m v_h \|_{L^{\infty}(e)}^2 \\
    & \leq C h_e^{d - 4} ( \| \mR_{K_1}^m v_h - v_{K_1} \|^2_{
    L^{\infty}(e)} + \| \mR_{K_2}^m v_h - v_{K_2} \|^2_{L^{\infty}(e)}
    ) \\
    & \leq C ( \| D^2(\mR_{K_1}^m v_h) \|_{L^2(K_1)}^2 + \|
    D^2(\mR_{K_2}^m v_h) \|_{L^2(K_2)}^2).
  \end{align*}
  For the boundary face $e \in \MEhB$, let $K$ be the element with $e
  \subset \partial K$ and we let $v_K := v_h|_K$. 
  Analogously, we deduce that 
  \begin{align*}
    h_e^{-1} \| \nabla_{\un} \mR^m v_h \|_{L^2(e)}^2  &\leq C
    h_e^{-1} ( \| \nabla_{\un} (\mR^m v_h - v_K) \|_{L^2(e)}^2 + \|
    \nabla_{\un} v_K \|_{L^2(e)}^2) \\
    & \leq C ( \| D^2 (\mR_K^m v_h) \|_{L^2(K)}^2 + h_e^{-1} \|
    \nabla_{\un} v_h \|_{L^2(e)}^2).
  \end{align*}
  Because $v_h$ vanishes on the boundary $\partial \Omega$, we have
  that
  \begin{align*}
    h_e^{-3} \| \mR^m v_h \|_{L^2(e)}^2 & =  h_e^{-3} \| \mR^m v_h -
    v_h\|_{L^2(e)}^2 \leq C h_e^{-3} \| \mR^m v_h - I_K (\mR^m v_h)
    \|_{L^2(e)}^2 \leq C \| D^2 (\mR^m_K v_h) \|_{L^2(K)}^2.
  \end{align*}
  Combining above estimates gives that 
  \begin{displaymath}
    \DGenorm{\mR^m v_h}^2 \leq C \sum_{K \in \MTh} \| D^2 (\mR^m
    v_h ) \|_{L^2(K)}^2 + C \DGenorm{v_h}^2.
  \end{displaymath}
  It remains to bound the $L^2$ norms for the second-order
  derivatives. 
  Notice that for any $K \in \MTh$, $v_h$ is piecewise
  linear in the domain $\mD(K)$.
  We let $K_{\max}, K_{\min} \in S(K)$ be elements such that 
  $| \nabla v_h|_{K_{\max}} - \nabla v_h|_{K_{\min}}|_{l^2} =
  \max_{K', K'' \in S(K)} | \nabla v_h|_{K'} - \nabla
  v_h|_{K''}|_{l^2}$.
  Let $v_{K, \max} := v_h|_{K_{\max}}$ and $v_{K, \min} :=
  v_h|_{K_{\min}}$ and we extend both linear polynomials $v_{K,
  \max}, v_{K, \min}$ to the domain $\mD(K)$ by the direct extension.
  Let $\tau_{\min}:= \min_{\bm{x} \in
  \mD(K)} (v_h - v_{K, \min} )$ and $\tau_{\max} := \max_{\bm{x} \in
  \mD(K)} (v_h - v_{K, \min})$,
  together with the definition \eqref{eq_LambdamK} and the inverse
  estimate, we obtain that
  \begin{align*}
    \| D^2 (\mR_K^m v_h) &\|_{L^2(K)}^2  = \| D^2(\mR^m v_h - v_{K,
    \min} - \tau_{\min}) \|_{L^2(K)}^2  \leq C h_K^{-4} \| \mR^m v_h -
    v_{K, \min} - \tau_{\min} \|_{L^{2}(K)}^2 \\
    & =  C h_K^{-4} \| \mR^m (v_h - v_{K, \min} - \tau_{\min})
    \|_{L^{2}(K)}^2 \leq C h_K^{d-4} \Lambda_{m, K}^2 \sum_{\bm{x} \in
    \mI(K)} ( v_h - v_{K, \min} - \tau_{\min} )^2 \\
    & \leq C  h_K^{d-4} \Lambda_{m, K}^2 \# \mI(K) (\tau_{\max} -
    \tau_{\min})^2.
  \end{align*}
  It is noted that $v_h - v_{h, \min}$ is also piecewise linear in
  $\mD(K)$, then $\tau_{\max}$ and $\tau_{\min}$ will be reached at
  two nodes $\bm{\nu}_0, \bm{\nu}_1$ in $\overline{\mD(K)}$. 
  We let $K^+, K^- \in \mS(K)$ be elements with $\bm{\nu}_0 \in
  \mN_{K^+}$ and $\bm{\nu}_1 \in \mN_{K^-}$. 
  There exist elements $K_1, \ldots, K_L(L = O(t_K))$ such that $K^+ =
  K_1, K^- = K_L$ and $K_j$ is face-adjacent to $K_{j+1}$, and we have
  that 
  \begin{align*}
    (\tau_{\max} - \tau_{\min})^2 \leq C h_K^2 \sum_{j = 1}^L \|
    \nabla (v_h - v_{K, \min}) \|_{L^{\infty}(K_j)}^2 \leq C t_K h_K^2
    \| \nabla v_{K_{\max}} - \nabla v_{K_{\min}} \|_{l^2}^2.
  \end{align*}
  For $K_{\max}, K_{\min}$, there exist elements $T_1,
  \ldots, T_J(J = O(t_K))$ such that $T_1 = K_{\min}, T_J = K_{\max}$
  and $T_j$ shares a common face $e_j$ with $T_{j+1}$.  Since $\nabla
  v_h$ is piecewise constant in $\mD(K)$ and $v_h$ is
  continuous on all faces, we deduce that 
  \begin{align*}
    \| \nabla v_h|_{K_1} - \nabla v_h|_{K_L} \|_{l^2}^2 \leq Ct_K
    \sum_{j = 1}^{L - 1} \| \nabla v_h|_{K_j} - \nabla v_h|_{K_{j+1}}
    \|_{l^2}^2 \leq  Ct_K \sum_{j = 1}^{L - 1} h_{e_j}^{-1} \|
    \jump{\nabla_{\un} v_h} \|_{L^2(e_j)}^2.
  \end{align*}
  Collecting all above estimates yields that 
  \begin{displaymath}
    \sum_{K \in \MTh} \| D^2(\mR^m_K v_h) \|_{L^2(K)}^2 \leq C
    \Lambda_m^2 \sum_{e \in \MEh} h_e^{-1} \| \jump{\nabla_{\un} v_h}
    \|_{L^2(e)}^2,
  \end{displaymath}
  which gives the desired estimate \eqref{eq_DGenorm_upper} and
  completes the proof.
\end{proof}

The condition number of the preconditioned system $A_{\mc{L}}^{-1}
A_m$ is estimated as below.
\begin{theorem}
  Let the penalty parameters $\mu_1$ and $\mu_2$ be taken as in Lemma
  \ref{le_ahbc}, then there exists a constant $C$ such that 
  \begin{equation}
    \kappa(A_{\mc{L}}^{-1} A_m) \leq C \Lambda_m^2.
    \label{eq_A1Amcond}
  \end{equation}
  \label{th_A1Amcond}
\end{theorem}
\begin{proof}
  In this proof, we associate any vector $\bmr{v} = \{ v_{\bm{\nu}}
  \}_{\bm{\nu} \in \MNhI} \in \mb{R}^{n_p}$ with a finite element
  function $v_h \in V_{h, 0}$ such that $v_h(\bm{\nu}) = v_{\bm{\nu}}$
  for $\forall \bm{\nu} \in \MNhI$.
  Now, we have that 
  \begin{displaymath}
    \bmr{v}^T A_{\mc{L}} \bmr{v} = a_h^{\mc{L}}(v_h, v_h), \quad
    \bmr{v}^T A_m \bmr{v} = a_h(\mR^m v_h, \mR^m v_h), \quad
    \forall \bmr{v} \in \mb{R}^{n_p}.
  \end{displaymath}
  By Lemma \ref{le_DGenorm_upper}, Lemma \ref{le_DGenorm_lower} and the 
  boundedness and the coercivity of $a_{h}(\cdot, \cdot)$, we find that 
  \begin{align*}
    {\Lambda}_m^{-2} \bmr{v}^T A_{m} \bmr{v} \leq C \Lambda_m^{-2}
    \DGenorm{\mR^m v_h}^2 \leq C \DGenorm{v_h}^2 = C a_h^{\mc{L}}(v_h,
    v_h) \leq C \DGenorm{\mR^m v_h}^2 \leq C  \bmr{v}^T A_{m} \bmr{v},
    \quad \forall \bmr{v} \in \mb{R}^{n_p}.
  \end{align*}
  From \cite[Lemma 2.1]{Xu1992iterative}, there holds
  $\kappa(A_{\mc{L}}^{-1} A_m) \leq C  \Lambda_m^2$, which completes
  the proof.
\end{proof}

From Theorem \ref{th_A1Amcond}, the resulting linear system $A_m
\bmr{x} = \bmr{b}$ can be solved by the CG method or other
Krylov-subspace family method (e.g. BiCGSTAB, GMRES), together with the
preconditioner $A_{\mc{L}}^{-1}$.
In the iteration, we are required to compute the matrix-vector product
$A_{\mc{L}}^{-1} \bmr{z}$, which is usually implemented by solving the
linear system $A_{\mc{L}} \bmr{y} = \bmr{z}$.
For this goal, we present a $\mc{W}$-cycle multigrid method for the
matrix $A_{\mc{L}}$ (see Algorithm~\ref{alg_mgsolver}),
following from the idea in the smoothed aggregation method
\cite{Vanek2001convergence, Vanek1996algebraic}.
Assume that the mesh $\MTh$ is obtained by successively refining a
coarse mesh $\mc{T}_1$ for several times, i.e. there exist a sequence
of meshes $\mc{T}_1, \mc{T}_2, \ldots, \mc{T}_J$ such that
$\mc{T}_{l+1}$ is created by subdividing all of triangular
(tetrahedral) elements in $\mc{T}_l$, and here $\mc{T}_J = \MTh$.
Let $h_j$ be the mesh size of $\mc{T}_j$, there holds
$h_j/h_k = 2^{k - j}$ and $h_J = h$.
Let $V_l \in H_0^1(\Omega)$ be the piecewise linear polynomial
space on the mesh $\mc{T}_l$, and we have the embedding relation 
$V_1 \subset V_{2} \subset \ldots \subset V_{J}$ with $V_J = V_{h, 0}$.
Let $I^j: L^2(\Omega) \rightarrow V_j$ be the $L^2$ projection
operator into the space $V_j$, and we write $I_{j}^{k}: V_j
\rightarrow V_{k}$ be the $L^2$ projection operator from $V_j$ to
$V_k$. 
On the finest level $J$, we define a linear operator 
$A: V_J \rightarrow V_J$ by
\begin{displaymath}
  (A v_h, w_h)_{L^2(\Omega)} = a_{\mc{L}}(v_h, w_h), \quad \forall
  v_h, w_h \in V_J, 
\end{displaymath}
with the induced norm $ \| v_h \|_{A}^2 := (A v_h,v_h)_{L^2(\Omega)}$.
On each level $j < J$, we define the operator $A_j$ and the
symmetric prolongator smoother $S_j$ in a recursive manner, which read
\begin{equation}
  A_{j} := (S_{j+1} I_{j}^{j+1})^T A_{j+1} S_{j+1} I_{j}^{j+1}, \quad
  S_j := I - 2.9 (\lambda_j)^{-1} A_j + 2.15 (\lambda_j)^{-2} A_j^2,
  \quad A_J := A,
  \label{eq_Al}
\end{equation}
where $\lambda_j$ is a constant that bounds $\varrho(A_j)$, i.e.
$\lambda_j \geq \varrho(A_j)$.
For any linear operator, we let
$\sigma(\cdot)$ and $\varrho(\cdot)$ be its spectrum and its largest
eigenvalue, respectively. 
The form of $S_j$ in \eqref{eq_Al} allows us to establish an estimate
to $\varrho(A_j)$ and show $S_j$ is positive definite.
We further demonstrate that $\lambda_j$ can be selected as 
\begin{equation}
  \lambda_j = 16^{j - J} \lambda, 
  \label{eq_lambdal}
\end{equation}
with $\lambda$ an available upper bound for $\varrho(A)$. 
By the estimate \eqref{eq_kappaA0}, we know that $\lambda =
O(h^{-4})$.  
The operator $A_j$ can be rewritten into the form $A_j = Q_j^T A_J
Q_j$, 
where the operator $Q_j: V_j \rightarrow V_J$ is defined by
\begin{displaymath}
  Q_j = S_{J}I_{J-1}^{J} S_{J-1} I_{J-2}^{J-1} \ldots
  S_{j+1}I_{j}^{j+1}, \quad 1 \leq j \leq J - 1, \quad Q_J = I. 
\end{displaymath}
By \cite{Vanek2001convergence}, the convergence analysis to 
Algorithm~\ref{alg_mgsolver} can be established on the following
results.

\begin{algorithm}[t]
  \caption{$\mc{W}$-cycle Multigrid Solver, MGSolverI($\bmr{x}_j$,
  $\bmr{b}_j$, $j$)}
  \label{alg_mgsolver}
  \KwIn{the initial guess $\bmr{x}_j$, the right hand side
  $\bmr{b}_j$, the level $j$;}
  \KwOut{the solution $\bmr{x}_j$;}
  \If{$j = 1$}{
  $A_1 \bmr{x}_1 = \bmr{b}_1$ is solved by the direct method. 

  return $\bmr{x}_1$;
  }

  \If{$j > 1$}{
  pre-smoothing: apply Gauss-Seidel sweep on $A_j \bmr{x}_j =
  \bmr{b}_j$;  

  correction on coarse mesh: set $\bmr{y} = (S_j
  I_{j-1}^{j})^T(\bmr{b}_j - A_j \bmr{x}_j)$;

  let $\bmr{z}_1 = \bmr{0}$, and update $\bmr{z}_1$ by
  MGSolver($\bmr{z}_1$, $\bmr{y}$, $j - 1$);

  set $\bmr{z}_2 = \text{MGSolver($\bmr{z}_1$, $\bmr{y}$, $j - 1$)}$; 

  set $\bmr{x}_j = \bmr{x}_j + S_j I_{j - 1}^j \bmr{z}_2$;

  post-smoothing: apply Gauss-Seidel sweep on $A_j \bmr{x}_j =
  \bmr{b}_j$;

  return $\bmr{x}_j$;
  }
\end{algorithm}

We first give the weak approximation property of the space $V_j$. 
\begin{lemma}
  There exists a constant $C$ such that
  \begin{equation}
    \min_{w_h \in V_j} \| v_h -  w_h  \|_{L^2(\Omega)}^2 \leq C 16^{J
    - j} h^4 \| v_h \|_{A}^2, \quad \forall v_h \in V_J, \quad 1 \leq
    j \leq J.
    \label{eq_weakerapp}
  \end{equation}
  \label{le_weakerapp}
\end{lemma}
\begin{proof}
  By \cite[Lemma 3]{Neilan2019discrete} and
  \cite[Lemma 3.1]{Georgoulis2011posteriori}, for the piecewise linear
  polynomial function $v_h$, there exists $v_h^{\mr{c}} \in
  H^2(\Omega)$ such that
  \begin{displaymath}
    \sum_{K \in \MTh} h_K^{2q-4} | v_h - v_h^{\mr{c}} |_{H^q(K)}^2
    \leq C \sum_{e \in \MEhI} h_e^{-1} \| \jump{\nabla_{\un} v_h}
    \|_{L^2(e)}^2 \leq C \| v_h \|_{A}^2, \quad 0 \leq q \leq 2.
  \end{displaymath}
  Because $v_h$ is piecewise linear, we have that $ |
  v_h^{\mr{c}}|_{H^2(\Omega)} \leq C \| v_h \|_{A}$. 
  Let $w_h := I^j v_h^{\mr{c}}$ be the $L^2$ interpolant of
  $v_h^{\mr{c}}$ into the space $V_j$. 
  By the approximation property of $I^j$, we have that 
  \begin{align*}
    \| v_h - w_h \|_{L^2(\Omega)}^2 & \leq C( \| v_h - v_h^{\mr{c}} 
    \|_{L^2(\Omega)}^2 + \| v_h^{\mr{c}}  - w_h \|_{L^2(
    \Omega)}^2) \leq C (h_J^4 \| v_h \|_{A}^2 + h_j^4 \| v_h
    \|_{A}^2)  \\
    & \leq C (h_j/h_J)^4 h_J^{4} \| v_h \|_{A}^2 = C16^{J - j} h^4 \|
    v_h \|_{A}^2,
  \end{align*}
  which brings us the estimate \eqref{eq_weakerapp} and
  completes the proof.
\end{proof}
Next, we show the selection \eqref{eq_lambdal} is proper.
\begin{lemma}
  Let $\lambda_j$ be taken as \eqref{eq_lambdal}, there holds 
  \begin{equation}
    \varrho(A_{j}) \leq \varrho(S_{j+1}^T A_{j+1} S_{j+1}) \leq 16^{j
    - J} \lambda = \lambda_j, \quad 1 \leq j \leq J - 1. 
    \label{eq_Ajupper} 
  \end{equation}
  \label{le_Ajupper}
\end{lemma}
\begin{proof}
  We begin by proving the estimate \eqref{eq_Ajupper} for the finest
  level $j = J - 1$.
  In this case, $S_{j+1}$ has the form $S_J = I - 2.9(\lambda_J)^{-1}
  A_J + 2.15(\lambda_J)^{-2} A_J^2$. 
  By the definition of $A_j = A_{J - 1}$, we deduce that 
  \begin{align*}
    \varrho(A_{J - 1}) & = \max_{v_h \in V_{J - 1}} \frac{ (Q_{J - 1} 
    v_h, A_J (Q_{J - 1}  v_h) )_{L^2(\Omega)}  }{ (v_h,
    v_h)_{L^2(\Omega)} } = \max_{v_h \in V_{J - 1} } \frac{ ( I_{J -
    1}^J v_h,
    (S_J^T A_J S_J ) (I_{J - 1}^J v_h))_{L^2(\Omega)}}{  (v_h,
    v_h)_{L^2(\Omega)} } \\ 
    & \leq  \max_{v_h \in V_J} \frac{ (v_h,
    (S_J^T A_JS_J) v_h)_{L^2(\Omega)}}{  (v_h,
    v_h)_{L^2(\Omega)} } \leq \varrho(S_J^TA_JS_J). 
  \end{align*}
  Since $\lambda_J = \lambda \geq \varrho(A_J)$, we know that
  $\varrho(A_J / \lambda_J) \leq 1$. 
  By the direct calculation, we obtain that
  \begin{displaymath}
    S_J^T A_J S_J = (I - 2.9\lambda_J^{-1} A_J  + 2.15\lambda_J^{-2}
    A_J^2)^2 A_J,
  \end{displaymath}
  and
  \begin{align*}
    \varrho(S_J^T A_J S_J) = \lambda_J \max_{ t \in
    \sigma(\lambda_J^{-1} A_J)} (1 - 2.9t + 2.15t^2)^2t \leq \lambda_J
    \max_{t \in [0,1]} (1 - 2.9t + 2.15t^2)^2t  \leq 16^{-1} \lambda_J,
  \end{align*}
  which gives the estimate \eqref{eq_Ajupper} for $j = J - 1$ and
  indicates $\varrho(A_{J - 1}) \leq \lambda_{J - 1}$. 
  The last inequality follows from the fact that the maximum value of
  the function $f(t) = (1 - 2.9t + 2.15t^2)^2t$ on $[0, 1]$ is reached
  at $t_0 = 1$, that is $\max_{t \in [0, 1]} f(t) \leq f(t_0) =
  16^{-1}$.
  By the same procedure, we can conclude the estimate
  \eqref{eq_Ajupper} at level $j - 1$ from the result at $j$.
  This completes the proof.
\end{proof}
Combining the estimate \eqref{eq_Ajupper} and the definition of $S_j$, 
we have that 
\begin{equation}
  \varrho(S_j) = \max_{t \in \sigma(\lambda_j^{-1} A_j)} (1 - 2.9t +
  2.15t^2) \leq \max_{t \in [0, 1]}  (1 - 2.9t +
  2.15t^2) = 1.
  \label{eq_vSj}
\end{equation}
It can be easily checked that on $[0, 1]$ the maximum value of $g(t)
:=  1 - 2.9t + 2.15t^2$
is achieved at $t_1 = 0$ while the minimum value of $g(t)$ is obtain
at $t_2 = \frac{2.9}{4.3}$.
This fact immediately brings that $0 < g(t_2) \leq g(t)(0 \leq t \leq
1) \leq g(t_1) =
1$, which gives the last equality in \eqref{eq_vSj} and 
\begin{equation}
  \sigma(S_j) \subset (0, 1], \quad 1 \leq j \leq J.
  \label{eq_sigmaSj}
\end{equation}
We can know that $S_j$ and $A_j$ are both symmetric and positive
definite.
From $A_j$, we introduce an induced norm $\| v_h \|_{A_j} := (v_h, A_j
v_h)_{L^2(\Omega)}$ for $\forall v_h \in V_j$. 
We further derive that
\begin{align}
  \|v_h - S_j v_h \|_{L^2(\Omega)} & = \| (2.9 \lambda_j^{-1} A_j - 2.15
  \lambda_j^{-2} A_j^2)v_h \|_{L^2(\Omega)} \leq C ( \lambda_j^{-1}\|
  A_j v_h \|_{L^2(\Omega)} + \lambda_j^{-2} \| A_j^2 v_h
  \|_{L^2(\Omega)}) \nonumber \\
  & \leq C \lambda_j^{-1/2} \| v_h \|_{A_j} \leq
  \frac{C}{\sqrt{\varrho(A_j)}} \| v_h \|_{A_j}, \quad \forall v_h \in
  V_j. \label{eq_ISj}
\end{align}
Let $\wt{Q}_j := Q_jI_J^j$ on the space $V_J$, we have the following
results.
\begin{lemma}
  There exist constants $C_0, C_1$ such that 
  \begin{equation}
    \| \wt{Q}_j v_h \|_{A} \leq \wt{C}(j) \| v_h \|_A, \quad \forall
    v_h \in V_J, \quad 1 \leq j \leq J,
    \label{eq_Qlbound}
  \end{equation}
  where $\wt{C}(j) = C_0 + C_1j$, and there exist constants $C_2, C_3$
  such that 
  \begin{equation}
    \|(\wt{Q}_{j} - \wt{Q}_{j+1}) v_h \|_{L^2(\Omega)} \leq
    \frac{\wh{C}(j)}{ \sqrt{\varrho(A_{j+1})}} \| v_h \|_A, \quad \forall
    v_h \in V_J, \quad
    1 \leq j \leq J - 1,
    \label{eq_QlQl1diff}
  \end{equation}
  where $\wh{C}(j) = C_2 + C_3j$.
  \label{le_MGconv}
\end{lemma}
\begin{proof}
  For the estimate \eqref{eq_Qlbound}, we have that 
  \begin{align}
    \| \wt{Q}_{j} v_h \|_{A} & = \| Q_{j} I_J^{j} v_h \|_A = \|
    Q_{j+1} S_{j+1} I_{j}^{j+1} I_J^{j} v_h \|_A = \| S_{j+1}
    I_{j}^{j+1} 
    I_{J}^{j} v_h \|_{A_{j+1}} \label{eq_wtQjvh} \\
    & \leq \| S_{j+1} I_J^{j+1} v_h\|_{A_{j+1}} + \| S_{j+1}
    (I_J^{j+1} - I_{j}^{j+1} I_J^{j}) v_h \|_{A_{j+1}}. \nonumber
  \end{align}
  Since $S_{j+1}$ is $A_{j+1}$-symmetric and $\varrho(S_{j+1}) \leq
  1$, we have that $\| S_{j+1} I_J^{j+1} v_h \|_{A_{j+1}} \leq \|
  I_J^{j+1} v_h \|_{A_{j+1}} = \| Q_{j+1} I_J^{j+1} v_h \|_{A} =
  \|\wt{Q}_{j+1} v_h \|_{A}$.
  We apply the approximation property \eqref{eq_weakerapp} and
  \eqref{eq_Ajupper} to find that
  \begin{align}
    \| S_{j+1}(I_J^{j+1} - &I_{j}^{j+1} I_{J}^{j}) v_h \|_{A_{j+1}}^2
    = ( (I_J^{j+1} - I_{j}^{j+1} I_{J}^{j}) v_h , S_{j+1}^T A_{j+1}
    S_{j+1} (I_J^{j+1} - I_{j}^{j+1} I_{J}^{j}) v_h )_{L^2(\Omega)}
    \nonumber \\
    & \leq 16^{j - J} \lambda \|  (I_J^{j+1} - I_{j}^{j+1} I_{J}^{j})
    v_h  \|_{L^2(\Omega)}^2
    \leq  16^{j - J} \lambda \|v_h - I_J^{j} v_h \|_{L^2(\Omega)}^2
    \leq C \| v_h \|_{A}^2.\label{eq_Sj1}
  \end{align}
  Combining \eqref{eq_wtQjvh} and \eqref{eq_Sj1} directly yields the
  estimate \eqref{eq_Qlbound}.

  We then focus on the second estimate \eqref{eq_QlQl1diff}.  
  From \eqref{eq_sigmaSj} and \eqref{eq_weakerapp}, we deduce that 
  \begin{align*}
    \| (\wt{Q}_{j} - \wt{Q}_{j+1}) v_h \|_{L^2(\Omega)} & = \| Q_{j+1}
    (S_{j+1} I_{j}^{j+1} I_{J}^{j} - I_J^{j+1}) v_h \|_{L^2(\Omega)}
    \leq \| (S_{j+1} I_{j}^{j+1} I_{J}^{j} - I_J^{j+1}) v_h
    \|_{L^2(\Omega)} \\
    & \leq \|
    S_{j+1} (I_j^{j+1} I_J^{j} - I_J^{j+1}) v_h
    \|_{L^2(\Omega)} + \|  (I_J^{j+1}  - S_{j+1} I_J^{j+1} ) v_h
    \|_{L^2(\Omega)}.
  \end{align*}
  The first term can be bounded by the weak approximation property
  \eqref{eq_weakerapp} and \eqref{eq_sigmaSj}, which reads
  \begin{align*}
    \| S_{j+1} (I_j^{j+1} I_J^{j} -& I_J^{j+1}) v_h \|_{L^2(\Omega)} 
    \leq \|  I_j^{j+1} (I_J^{j} v_h -  v_h )\|_{L^2(\Omega)} \leq
    \| v_h - I_J^j v_h \|_{L^2(\Omega)}  \\
    & \leq C 16^{j - J} \lambda \| v_h \|_{A} \leq \frac{C}{
    \sqrt{\varrho(A_{j+1})} } \| v_h \|_{A}.
  \end{align*}
  The second term can be estimated by \eqref{eq_ISj}, and we find that
  \begin{align*}
    \| I_J^{j+1} v_h - S_{j+1} I_{J}^{j+1} v_h \|_{L^2(\Omega)} & \leq
    \frac{C}{\sqrt{\varrho(A_{j+1})}} \| I_J^{j+1} v_h \|_{A_{j+1}} =
    \frac{C}{\sqrt{\varrho(A_{j+1})}} \| Q_{j+1} I_J^{j+1} v_h \|_{A}
    \\
    & =  \frac{C}{\sqrt{\varrho(A_{j+1})}} \| \wt{Q}_{j+1} v_h \|_{A}
    \leq  \frac{C \wt{C}(j)}{\sqrt{\varrho(A_{j+1})}} \|  v_h \|_{A}.
  \end{align*}
  Collecting the above two estimates leads to the estimate
  \eqref{eq_QlQl1diff}, which completes the proof.
\end{proof}

By \cite[Theorem 3.5]{Vanek2001convergence}, Lemma \ref{le_MGconv}
gives the convergence rate of Algorithm~\ref{alg_mgsolver}, which
reads
\begin{equation}
  \| \bmr{x} - \text{MG}(\bmr{x}, \bmr{b})\|_{A} \leq ( 1-
  \frac{1}{C(J)})
  \| \bmr{x} - \bmr{x} \|_{A}, \quad 
  \label{eq_MGconv}
\end{equation}
where $C(J) = O(J^3)$.

Consequently, the linear system $A_m \bmr{x} = \bmr{b}$ of
\eqref{eq_dvar} can be solved by the CG method using
Algorithm~\ref{alg_mgsolver} as the preconditioner.

Ultimately, we present another $\mc{W}$-cycle
multigrid algorithm for the system $A_{\mc{L}} \bmr{y} = \bmr{z}$,
which has a simpler implementation and is more close to the
geometrical multigrid method, see Algorithm~\ref{alg_mgsolverII}.
The major difference with Algorithm~\ref{alg_mgsolver} is the smoother
matrix $S_j$ is not required.
Algorithm~\ref{alg_mgsolverII} is also based on the spaces $V_1
\subset V_2 \subset \ldots \subset V_J$. 
For every $\mc{T}_j$, we let $\mc{E}_j$ as the set of all $d-1$
dimensional faces in $\mc{T}_j$. 
On each level $j$, we define a bilinear form 
\begin{displaymath}
  a_{\mc{L}, j}(v_h, w_h) := \sum_{e \in \mE_j} \int_e h_e^{-1}
  \jump{\nabla_{\un} v_h} \jump{\nabla_{\un} w_h} \dx{s}, \quad
  \forall v_h, w_h \in V_j,
\end{displaymath}
and we have that $a_{\mc{L}, J}(\cdot, \cdot) = a_{\mc{L}}(\cdot,
\cdot)$.
This form $a_{\mc{L}, j}(\cdot, \cdot)$ comes from the interior
penalty scheme on the mesh $\mc{T}_j$ over the piecewise linear spaces
$V_j \times V_j$.
Let $A_{\mc{L}, j}$ be the matrix from $a_{\mc{L}, j}(\cdot, \cdot)$.
In Algorithm~\ref{alg_mgsolverII}, 
the prolongation operator and the restriction operator are standard,
and on each level we consider the linear system of the matrix
$A_{\mc{L}, j}$.
We note that $a_{\mc{L}, j}(\cdot, \cdot)$ has the form
\begin{displaymath}
  a_{\mc{L}, j}(v_h, w_h) = 2^{j - J} a_{\mc{L}}(v_h, w_h), \quad
  \forall v_h, w_h \in V_j, \quad 1 \leq j \leq J.
\end{displaymath}
Our numerical observation demonstrates that
Algorithm~\ref{alg_mgsolverII} also works well as the preconditioner.
The convergence study of Algorithm~\ref{alg_mgsolverII} is now left as
a future research.

\begin{algorithm}[t]
  \caption{$\mc{W}$-cycle Multigrid Solver, MGSolverII($\bmr{x}_j$,
  $\bmr{b}_j$, $j$)}
  \label{alg_mgsolverII}
  \KwIn{the initial guess $\bmr{x}_j$, the right hand side
  $\bmr{b}_j$, the level $j$;}
  \KwOut{the solution $\bmr{x}_j$;}
  \If{$j = 1$}{
  $A_{\mc{L}, 1} \bmr{x}_1 = \bmr{b}_1$ is solved by the direct method. 

  return $\bmr{x}_1$;
  }

  \If{$j > 1$}{
  pre-smoothing: apply Gauss-Seidel sweep on $A_{\mc{L}, j} \bmr{x}_j =
  \bmr{b}_j$;  

  correction on coarse mesh: set $\bmr{y} = 
  I_{j}^{j-1}(\bmr{b}_j - A_{\mc{L}, j} \bmr{x}_j)$;

  let $\bmr{z}_1 = \bmr{0}$, and update $\bmr{z}_1$ by
  MGSolverII($\bmr{z}_1$, $\bmr{y}$, $j - 1$);

  set $\bmr{z}_2 = \text{MGSolverII($\bmr{z}_1$, $\bmr{y}$, $j - 1$)}$; 

  set $\bmr{x}_j = \bmr{x}_j + I_{j - 1}^j \bmr{z}_2$;

  post-smoothing: apply Gauss-Seidel sweep on $A_{\mc{L}, j} \bmr{x}_j
  = \bmr{b}_j$;

  return $\bmr{x}_j$;
  }
\end{algorithm}

\section{Numerical Results}
\label{sec_numericalresults}
\begin{figure}[htb]
  \centering
  \includegraphics[width=0.23\textwidth]{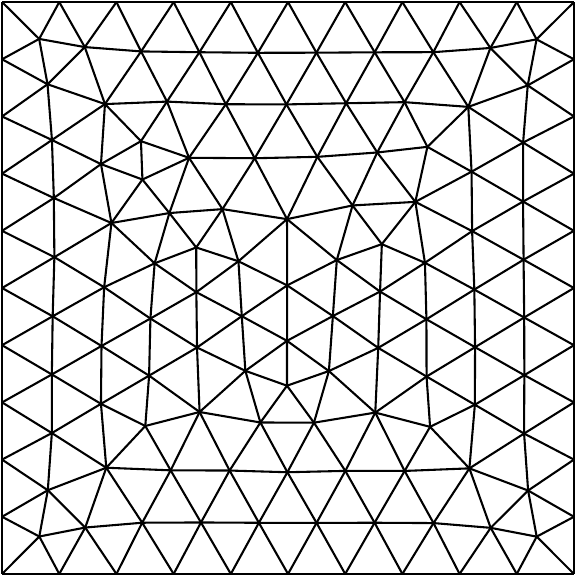}
  \hspace{100pt}
  \includegraphics[width=0.21\textwidth]{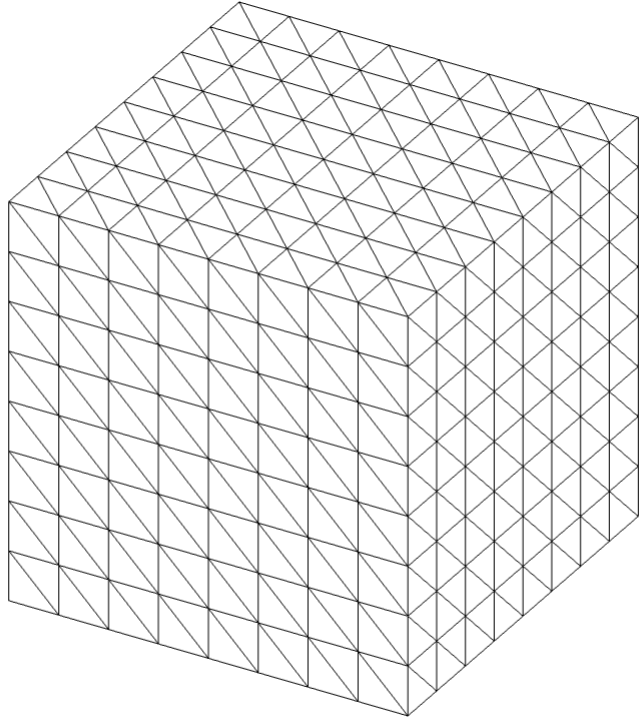}
  \caption{The triangular mesh (left) and tetrahedral mesh (right).}
  \label{fig_mesh}
\end{figure}

In this section, numerical experiments in both two and three
dimensions are conducted to demonstrate the performance of the
numerical scheme and the preconditioning method.
We adopt a family of quasi-uniform triangular and tetrahedral meshes
to solve the problems in two and three dimensions, respectively, see
Fig.~\ref{fig_mesh}.

\begin{figure}[htb]
  \centering
  \includegraphics[width=0.25\textwidth]{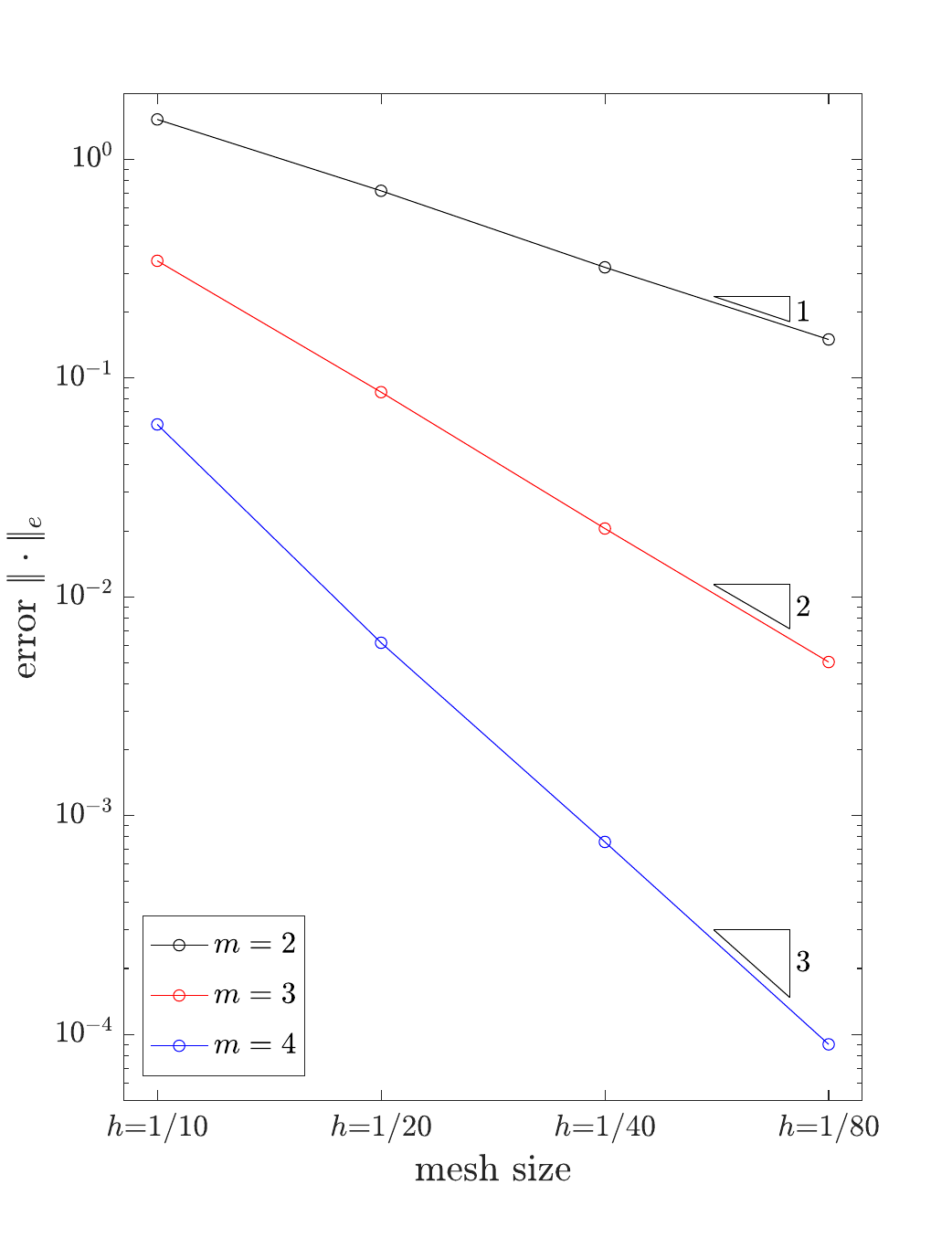}
  \hspace{50pt}
  \includegraphics[width=0.25\textwidth]{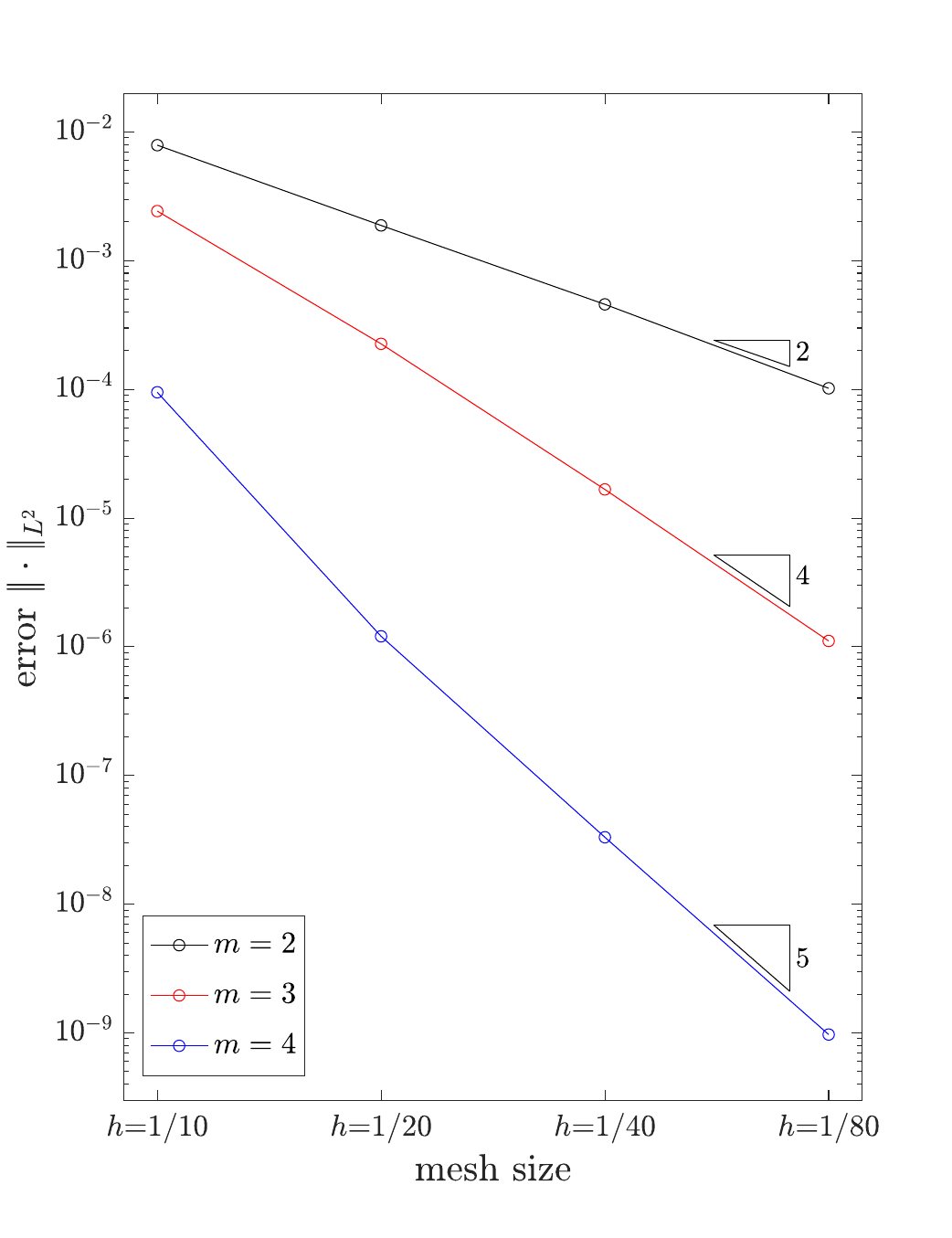}
  \caption{The numerical errors under the energy norm (left)/$L^2$
  norm (right) in Example 1.}
  \label{fig_ex1err}
\end{figure}

\noindent\textbf{Example 1.}
In the first example, we solve a biharmonic equation in two dimensions
defined on the squared domain $\Omega = (0, 1)^2$. 
The exact solution is given by 
\begin{displaymath}
  u(x,y) = \sin^2(\pi x) \sin^2(\pi y).
\end{displaymath}
We apply the approximation space $U_h^m$ of the accuracy $m = 2, 3, 4$
to solve the biharmonic problem on the triangular meshes with the mesh
size $h = 1/10, 1/20, 1/40, 1/80$.
We first test the accuracy of the scheme.  
The numerical errors under both the energy norm and the $L^2$ norm 
are plotted in Fig.~\ref{fig_ex1err}. 
It can be seen that the numerically detected convergence rates
validate the theoretical preconditions in Theorem \ref{th_DGerror}.

\begin{figure}[htb]
  \centering
  \includegraphics[width=0.25\textwidth]{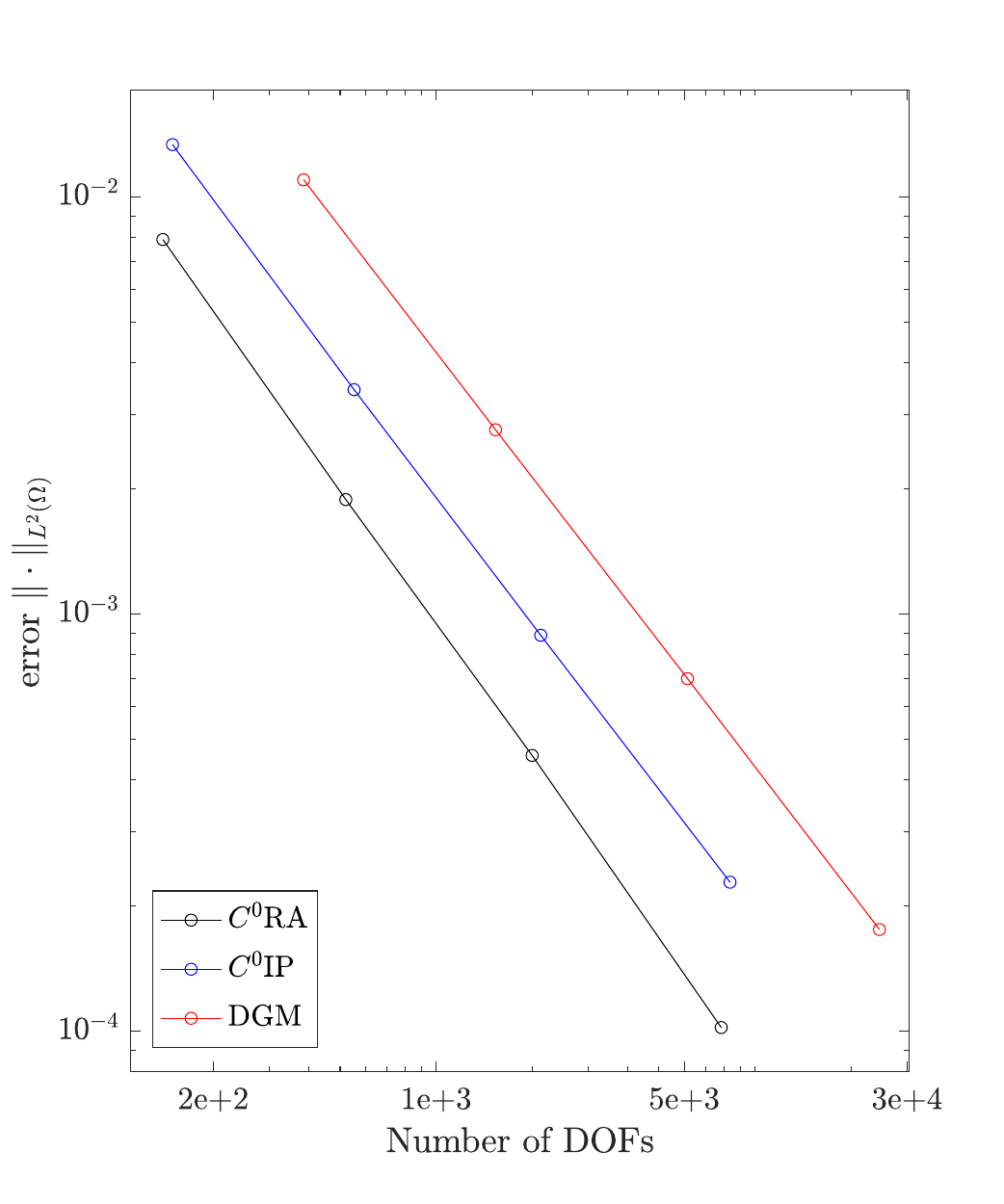}
  \hspace{20pt}
  \includegraphics[width=0.25\textwidth]{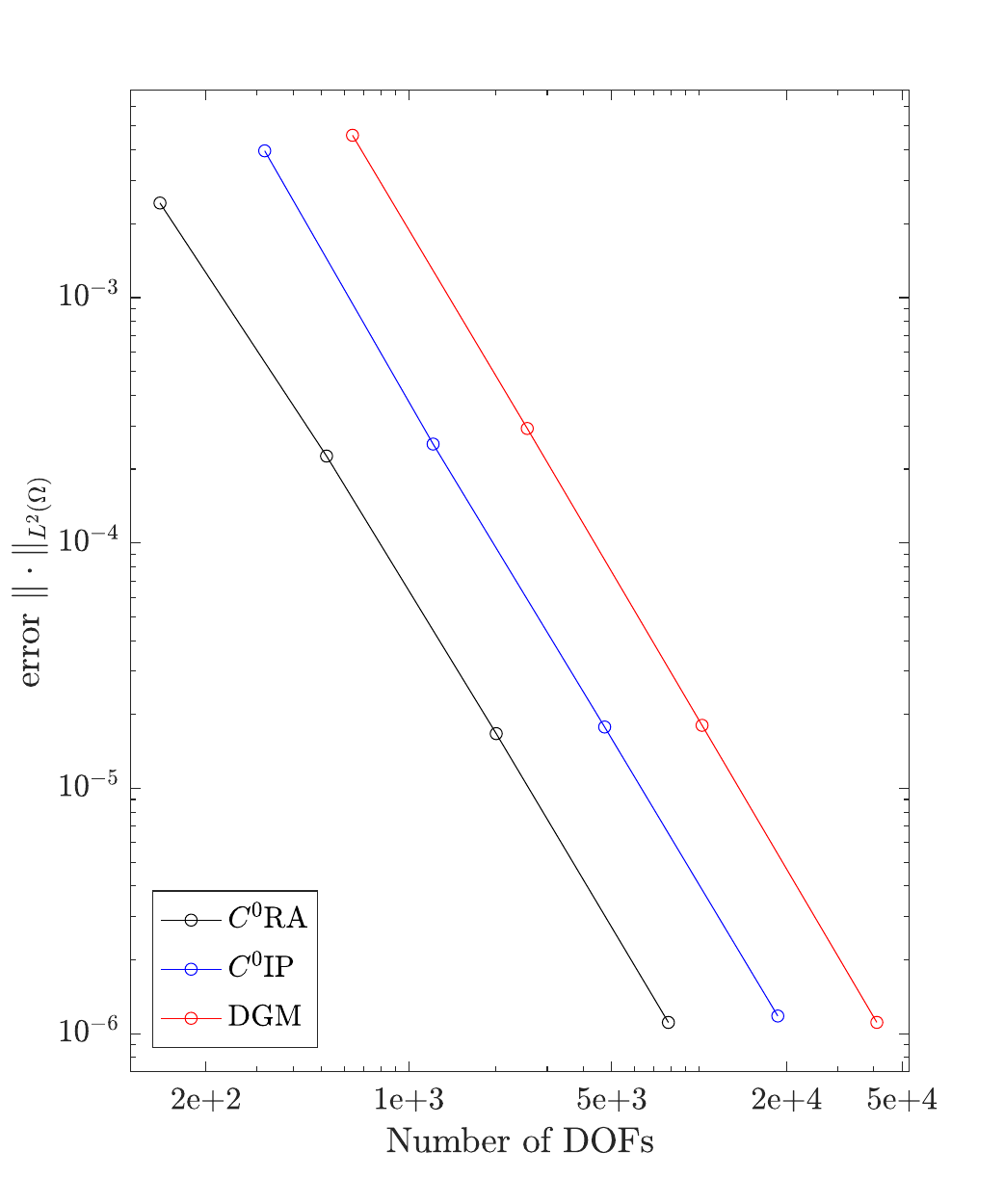}
  \hspace{20pt}
  \includegraphics[width=0.25\textwidth]{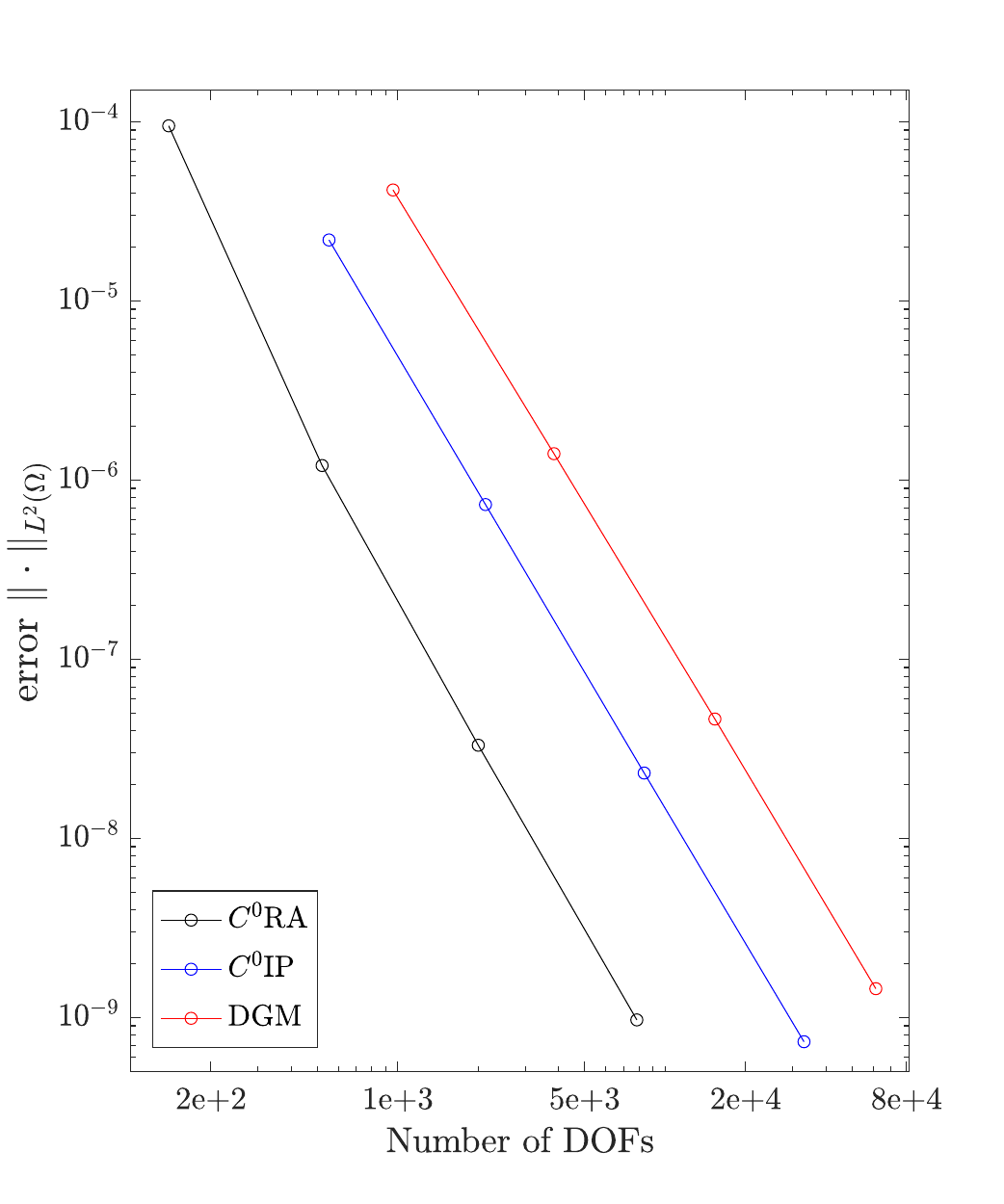}
  \caption{The numerical comparison of $L^2$ errors for different
  methods (from left to right, $m = 2, 3, 4$).}
  \label{fig_ex1comp}
\end{figure}

We also numerically compare the proposed method with the interior
penalty discontinuous Galerkin method and the $C^0$ interior penalty
method. 
In Fig.~\ref{fig_ex1comp}, we plot the $L^2$ errors of solutions from
different methods against the number of degrees of freedom.
It can be seen that all methods have same convergence rates, which are
consistent with the error estimates. 
To achieve the same $L^2$ error, the reconstructed method uses less 
degrees of freedom and 
is observed to be more computationally efficient than other methods.

We further present the numerical performance of the preconditioning
method in solving the resulting linear system. 
The condition numbers of the system $A_m$ and the
preconditioned system $A_{\mc{L}}^{-1} A_m$ are listed in
Tab.~\ref{tab_cnexample1}.
It is clear that $\kappa(A_m)$ grows very fast at
the speed of $O(h^{-4})$ while for all $m$, $\kappa(A_{\mc{L}}^{-1}
A_m)$ increases very slightly as the mesh size tends to zero.  
The numerical observation illustrates the results given in Theorem
\ref{th_A1Amcond}.
As stated in the previous section, we apply the preconditioned CG
method using the approximation of $A_{\mc{L}}^{-1}$ as the
preconditioner, i.e. Algorithm \ref{alg_mgsolver} and Algorithm
\ref{alg_mgsolverII}, for the resulting linear system.
The stopping criteria is set by $\| \bmr{r}_k \|_{l^2} \leq 10^{-9} \|
\bmr{b} \|_{l^2}$, where $\bmr{r}_k := \bmr{b} - A_m \bm{x}_k$ is the
residual at step $k$.
The iteration counts for PCG/CG methods are gathered in
Tab.~\ref{tab_ex1steps}.
For the standard CG method, the iteration count quadruples after the
mesh is refined, which agrees with the estimation to the condition
number. 
For the preconditioning methods, the required steps are numerically
observed to grow slightly when the mesh size approaches zero.
The numerical results illustrate the efficiency of our method.

\begin{table}[htb]
  \centering
  \renewcommand{\arraystretch}{1.25}
  \scalebox{0.78}{
  \begin{tabular}{p{0.5cm}|p{2.0cm}|p{1.8cm}|p{1.8cm}|p{1.8cm}| 
    p{1.8cm}}
    \hline\hline
    \multicolumn{2}{c|}{\diagbox[width=3.2cm]{$h$}{$m$}} 
    & 1/10 & 1/20 & 1/40 & 1/80 \\
    \hline
    \multirow{2}{*}{2} & $\kappa(A_{\mc{L}}^{-1} A_m)$ & 
    10.07 & 15.00 & 16.12 & 17.23 \\ 
    \cline{2-6}
    &$\kappa(A_m)$ & 3.673e+3 & 9.217e+4 & 1.530e+6 & 2.735e+7 \\ 
    \hline
    \multirow{2}{*}{3} & $\kappa(A_{\mc{L}}^{-1} A_m)$ &
    16.69 & 21.15 & 23.09 & 25.02 \\ 
    \cline{2-6}
    &$\kappa(A_m)$ 
    & 6.672e+3 & 1.256e+5 & 2.353e+6 & 4.173e+7 \\ 
    \hline
    \multirow{2}{*}{4} & $\kappa(A_{\mc{L}}^{-1} A_m)$ 
    & 32.75 & 49.83 & 60.37 & 65.45 \\ 
    \cline{2-6}
    &$\kappa(A_m)$ & 1.582e+4 & 2.662e+5 & 4.26e+6 & 7.172e+7 \\ 
    \hline \hline
    \end{tabular}
   }
  \caption{The condition numbers in Example 1.}
  \label{tab_cnexample1}
\end{table}

\begin{table}[htb]
  \centering
  \renewcommand{\arraystretch}{1.25}
  \scalebox{0.78}{
  \begin{tabular}{p{0.25cm}|p{3.5cm}|p{1.25cm}|p{1.25cm}|p{1.25cm}|
    p{1.25cm}}
    \hline\hline
    \multirow{2}{*}{$m$} & \diagbox[width=3.85cm]{Preconditioner}{$h$}
    & 1/10 & 1/20 & 1/40 & 1/80 \\ \hline
    \multirow{3}{*}{2} 
    & Algorithm \ref{alg_mgsolver} for $A_{\mc{L}}^{-1}$ 
    & 21 & 28 & 29 & 30 \\ 
    \cline{2-6}
    & Algorithm \ref{alg_mgsolverII}  for $A_{\mc{L}}^{-1}$ 
    & 21 & 27 & 29 & 31 \\ 
    \cline{2-6}
    & Identity  & 109 & 538 & 2338 & $>3000$\\ \hline        
    \multirow{3}{*}{3} 
    & Algorithm \ref{alg_mgsolver} for $A_{\mc{L}}^{-1}$ 
    & 31 & 45 & 51 & 54 \\ 
    \cline{2-6}
    & Algorithm \ref{alg_mgsolverII} for $A_{\mc{L}}^{-1}$
    & 31 & 46 & 49 & 51 \\ 
    \cline{2-6}
    & Identity 
    & 109 & 649 & $>3000$ & $ > 3000$ \\ \hline        
    \multirow{3}{*}{4} 
    & Algorithm \ref{alg_mgsolver} for $A_{\mc{L}}^{-1}$ 
    & 46 & 58 & 71 & 73 \\ 
    \cline{2-6}
    & Algorithm \ref{alg_mgsolverII} for $A_{\mc{L}}^{-1}$
    & 46  & 62  & 67 & 71  \\ 
    \cline{2-6}
    & Identity
    & 118 & 847 & $>3000$ & $>3000$ \\ 
    \hline\hline
  \end{tabular}
  }
  \caption{The convergence steps for PCG/CG methods in Example 1.}
  \label{tab_ex1steps}
\end{table}

\begin{figure}[htb]
  \centering
  \includegraphics[width=0.25\textwidth]{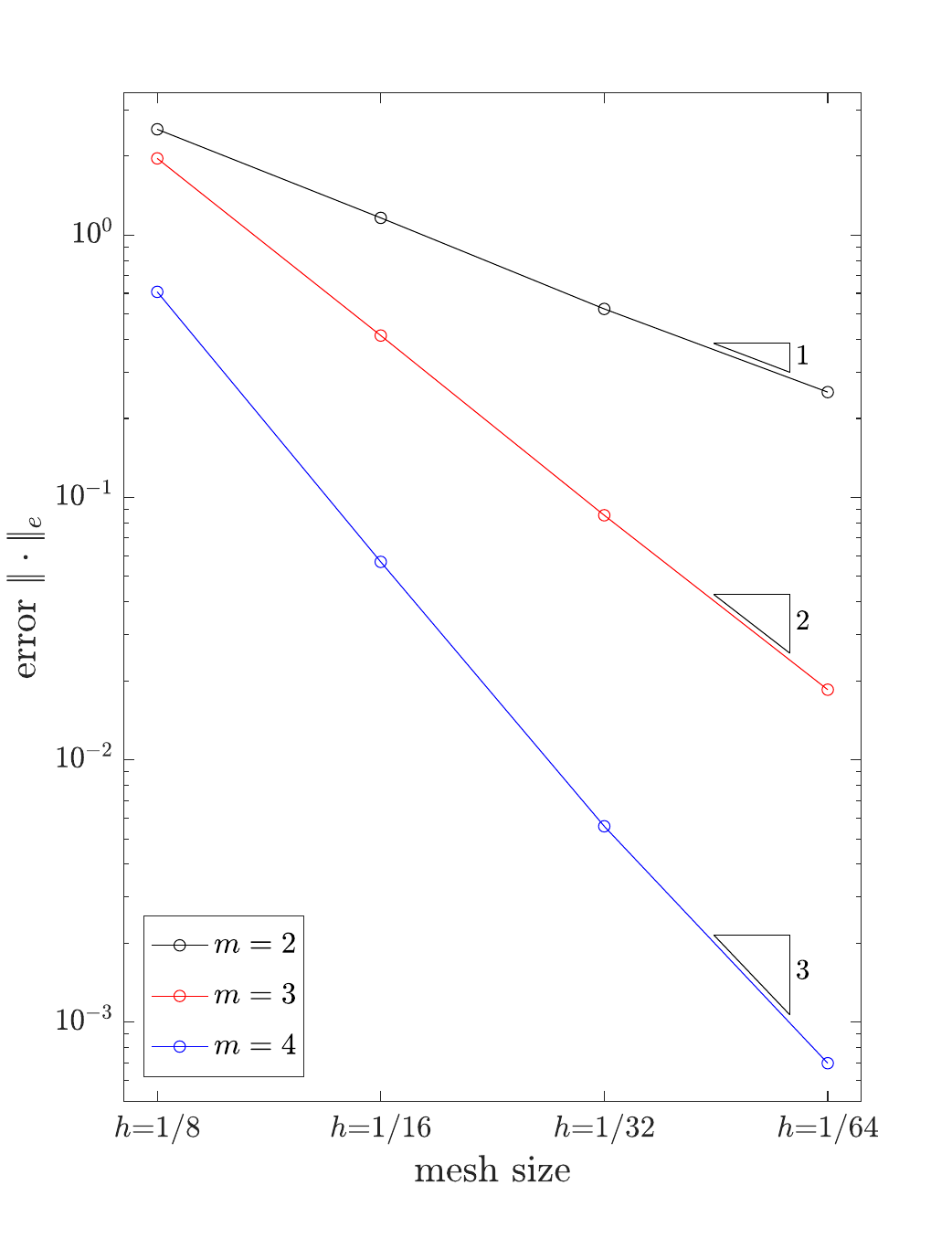}
  \hspace{50pt}
  \includegraphics[width=0.25\textwidth]{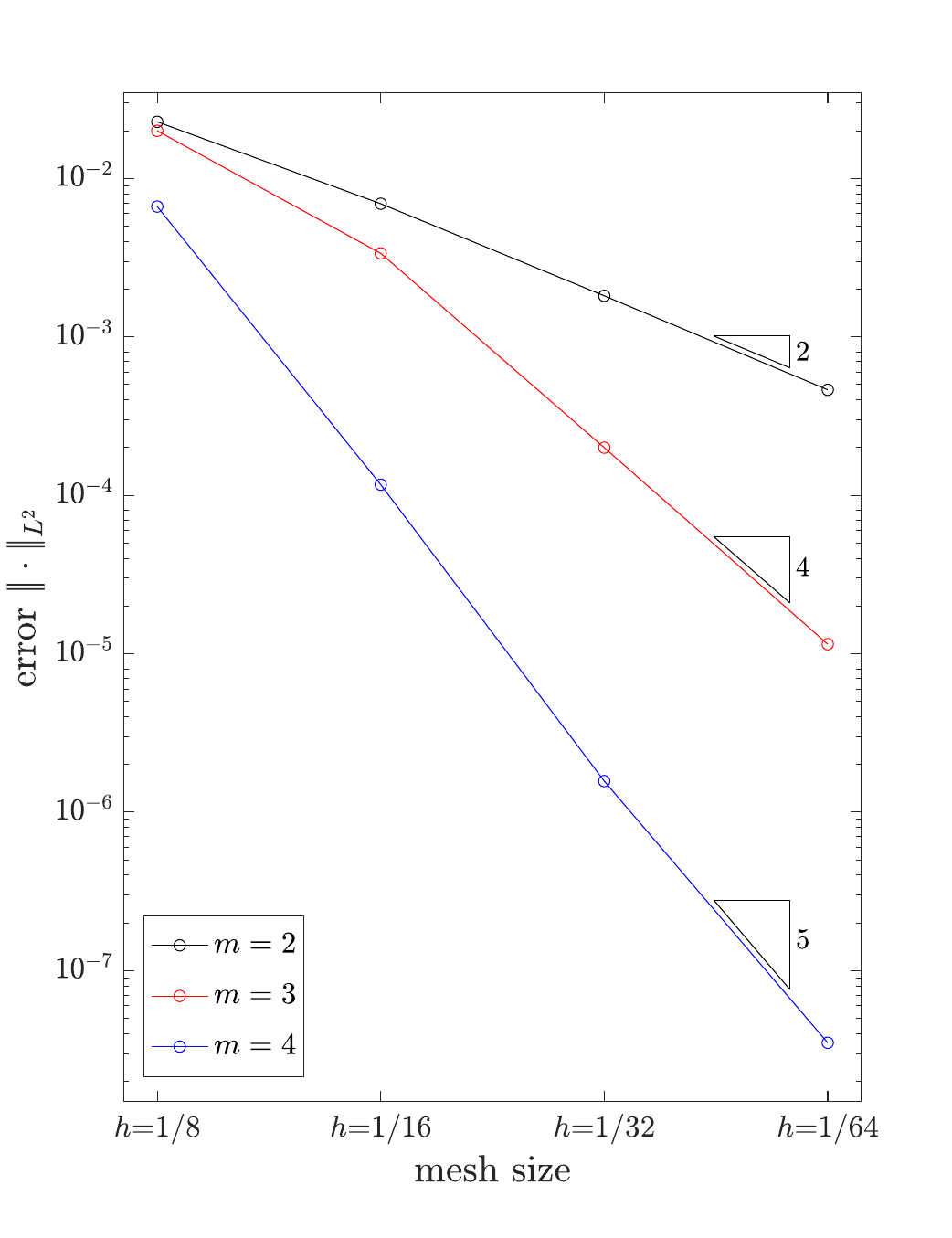}
  \caption{The numerical errors under the energy norm (left)/$L^2$
  norm (right) in Example 2.}
  \label{fig_ex2err}
\end{figure}

\noindent \textbf{Example 2.}
In this example, we solve a three-dimensional problem in the cubic
domain $\Omega = (0, 1)^3$.
The analytic solution is chosen to be 
\begin{displaymath}
  u(x,y,z) = \sin(\pi x) \sin(\pi y) \sin(\pi z).
\end{displaymath}
The source function and boundary conditions are taken accordingly.
This problem is discretized on a series of tetrahedral meshes with the
mesh size $h = 1/8, 1/16, 1/32, 1/64$. 
The numerical errors under $L^2$ norm and the energy norm are 
are shown in Fig.~\ref{fig_ex2err}.
As predicted, the numerically detected errors decrease sharply for the
proposed scheme, which agree with the theoretical estimation.

For the linear systems $A_m$ and $A_{\mc{L}}^{-1} A_m$, their
condition numbers are collected in Tab.~\ref{tab_cnexample2}.
Similar to the case of two dimensions, it can be seen that
$\kappa(A_m)$ grows at the rate $O(h^{-4})$ while
$\kappa(A_{\mc{L}}^{-1} A_m)$ is nearly constant. 
The convergence steps of linear solvers are gathered in
Tab.~\ref{tab_ex2steps}. 
In three dimensions, the linear solver can be
significantly accelerated with both multigrid methods for
$A_{\mc{L}}^{-1}$.
The numerical results again illustrate the efficiency of the
preconditioning method.

\begin{table}[htb]
  \centering
  \renewcommand{\arraystretch}{1.25}
  \scalebox{0.78}{
  \begin{tabular}{p{0.5cm}|p{2.0cm}|p{1.8cm}|p{1.8cm}|p{1.8cm}| 
    p{1.8cm}}
    \hline\hline
    \multicolumn{2}{c|}{\diagbox[width=3.85cm]{$h$}{$m$}} 
    & 1/8 & 1/16 & 1/32 & 1/64 \\ \hline
    \multirow{2}{*}{2} & $\kappa(A_{\mc{L}}^{-1} A_m)$ & 
    6.64 & 13.26 & 17.89 & 18.47 \\ 
    \cline{2-6}
    &$\kappa(A_m)$ & 1.013e+2 & 1.152e+3 & 1.957e+4 & 3.703e+5 \\ 
    \hline
    \multirow{2}{*}{3} & $\kappa(A_{\mc{L}}^{-1} A_m)$ &
    9.94 & 26.29 & 38.10 & 42.03 \\
    \cline{2-6}
    &$\kappa(A_m)$ & 2.286e+2 & 5.012e+3 & 1.106e+5 & 1.95e+6 \\ 
    \hline
    \multirow{2}{*}{4} & $\kappa(A_{\mc{L}}^{-1} A_m)$ & 
    40.50 & 158.47 & 219.59 & 227.29 \\ \cline{2-6}
    &$\kappa(A_m)$ & 
    5.164e+2 & 3.652e+4 & 6.292e+5 & 1.136e+7 \\ 
    \hline \hline
    \end{tabular}
   }
  \caption{The condition numbers in Example 2.}
  \label{tab_cnexample2}
\end{table}

\begin{table}[htb]
  \centering
  \renewcommand{\arraystretch}{1.25}
  \scalebox{0.78}{
  \begin{tabular}{p{0.25cm}|p{3.5cm}|p{1.25cm}|p{1.25cm}|p{1.25cm}|
    p{1.25cm}}
    \hline\hline
    \multirow{2}{*}{$m$} & \diagbox[width=3.85cm]{Preconditioner}{$h$} & 
    1/8 & 1/16 & 1/32 & 1/64 \\ \hline
    \multirow{3}{*}{2} 
    & Algorithm \ref{alg_mgsolver} for $A_{\mc{L}}^{-1}$
    & 25 & 38 & 48 & 52 \\ 
    \cline{2-6}
    & Algorithm \ref{alg_mgsolverII} for $A_{\mc{L}}^{-1}$ 
    & 25  & 37  & 47 & 49  \\ 
    \cline{2-6}
    & Identity & 99 & 366 & 1391 & $>3000$ \\ 
    \hline        
    \multirow{3}{*}{3}
    & Algorithm \ref{alg_mgsolver} for $A_{\mc{L}}^{-1}$
    & 29 & 50 & 63 & 69  \\ 
    \cline{2-6}
    & Algorithm \ref{alg_mgsolverII} for $A_{\mc{L}}^{-1}$ 
    & 28 & 49 & 65 & 71 \\ 
    \cline{2-6}
    & Identity & 189 & 996 & $> 3000$ & $>3000$ \\ \hline        
    \multirow{3}{*}{4} 
    & Algorithm \ref{alg_mgsolver} for $A_{\mc{L}}^{-1}$
    & 53  & 106 & 139 & 143  \\
    \cline{2-6}
    & Algorithm \ref{alg_mgsolverII} for $A_{\mc{L}}^{-1}$
    & 53 & 108  & 137  & 141  \\ 
    \cline{2-6}
    & Identity & 213 & 1094  & $>3000$ & $> 3000$ \\  
    \hline\hline
  \end{tabular}
  }
  \caption{The convergence steps for PCG/CG methods in Example 2.}
  \label{tab_ex2steps}
\end{table}

\section*{Acknowledgements}
This work was supported by National Natural Science Foundation of
China (no. 12201442).


\begin{appendix}
  \section{}
  \label{sec_reconM}
  In this appendix, we present the method for computing the constant
  $\Lambda_m$ and $\Lambda_{m, K}$ on each element after constructing
  element patches for a given mesh $\MTh$. 
  For every $K \in \MTh$, we let $p_1, p_2, \ldots, p_l$ be a group of
  standard orthogonal basis functions in $\mP_m(K)$ under the $L^2$
  inner product $(\cdot, \cdot)_{L^2(K)}$.
  Then, any $q \in \mP_m(K)$ can be expressed by a group coefficients
  $\bm{\alpha} = \{a_j\}_{j = 1}^l \in \mb{R}^{l}$ such that  $q =
  \sum_{j = 1}^l a_j p_j$.
  In addition, $q$ and all $p_j$ can be naturally extended to the
  domain $\mD(K)$.
  The main step to obtain $\Lambda_m$ is to compute $\Lambda_{m, K}$
  per element. 
  By $\{p_j \}_{j = 1}^l$, $\Lambda_{m, K}$ can be written into an
  algebraic form as
  \begin{displaymath}
    \Lambda_{m, K}^2 = \max_{\bm{\alpha} \in \mb{R}^l} \frac{|
    \bm{\alpha}|_{l^2}^2 }{h_K^d \bm{\alpha}^T B_K \bm{\alpha}},
    \quad B_K = \{b_{ij}\}_{l \times l}, \quad b_{ij} = \sum_{\bm{x}
    \in \mI(K)} p_i(\bm{x})p_j(\bm{x}).
  \end{displaymath}
  From the above matrix representation, we have that $\Lambda_{m, K} =
  (h_K^d \sigma_{\min}(B_K))^{-1/2}$, where $\sigma_{\min}(B_K)$ is
  the smallest singular value to the matrix $B_K$. 
  Therefore, it is enough to observe all the smallest singular values
  $\sigma_{\min}(B_K)$, and $\Lambda_m$ can be further obtained by
  \eqref{eq_LambdamK}.

\begin{figure}[htp]
    \centering
    \includegraphics[width=0.23\textwidth]{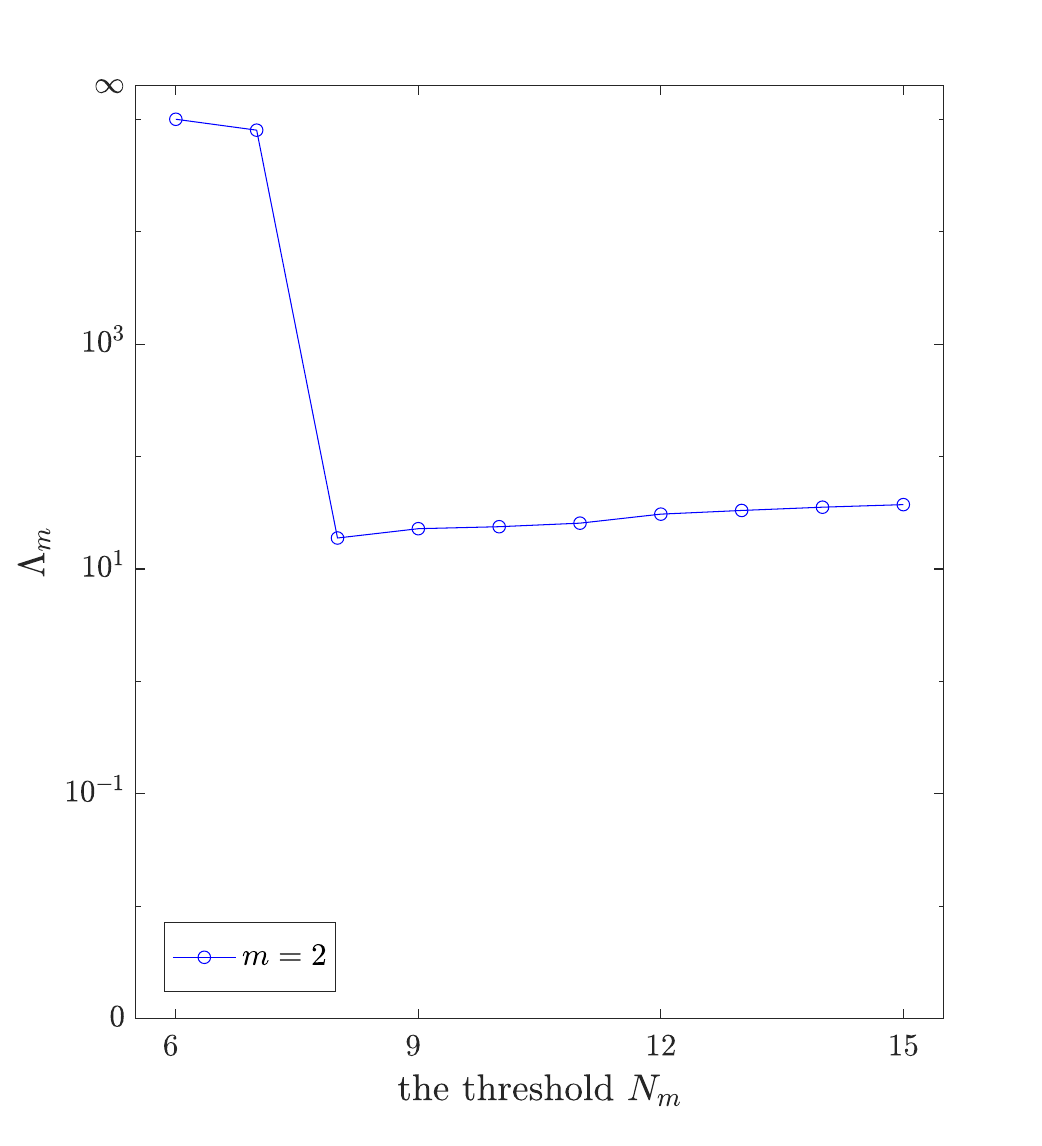}
    \hspace{10pt}
    \includegraphics[width=0.23\textwidth]{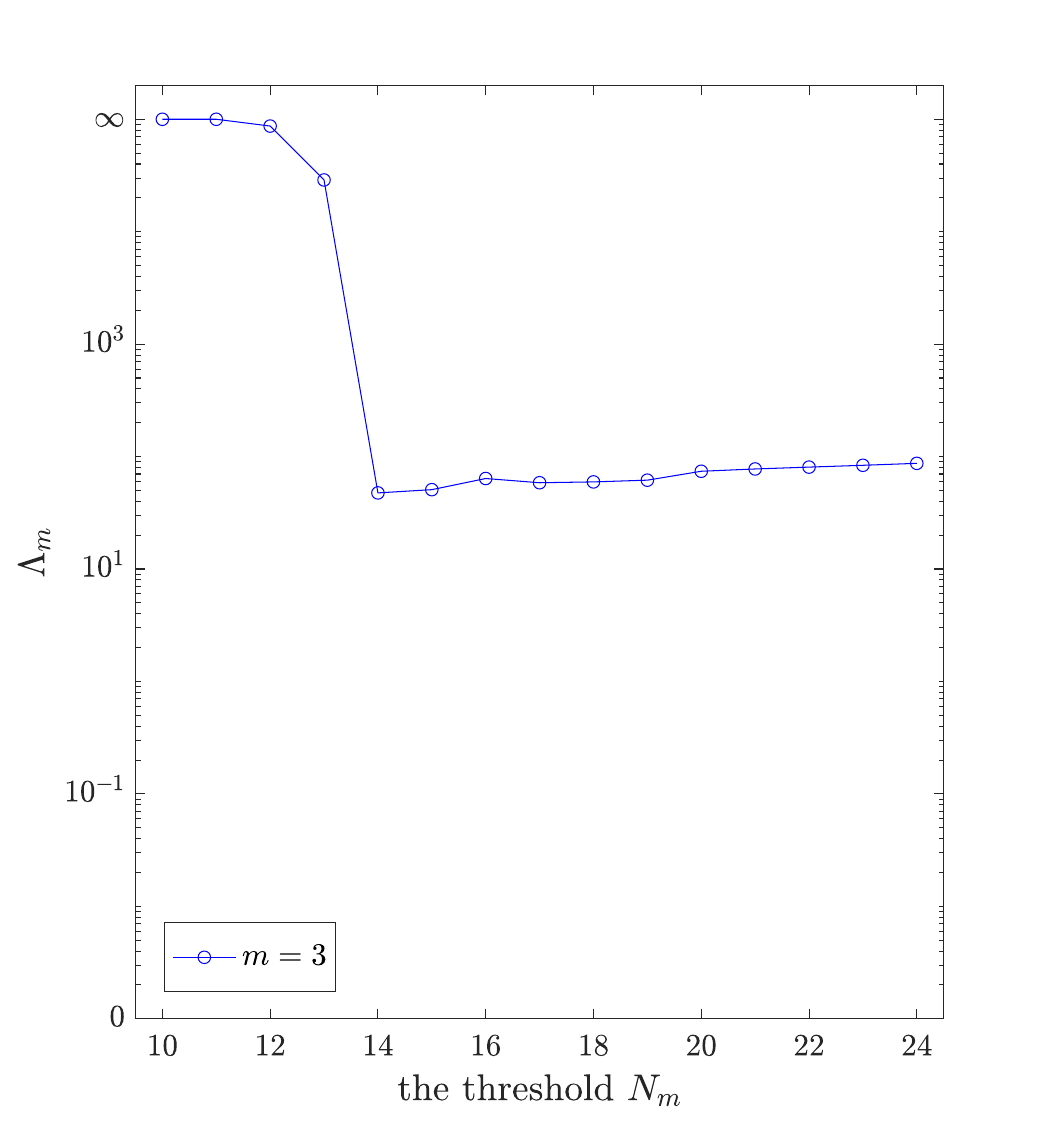}
    \hspace{10pt}
    \includegraphics[width=0.23\textwidth]{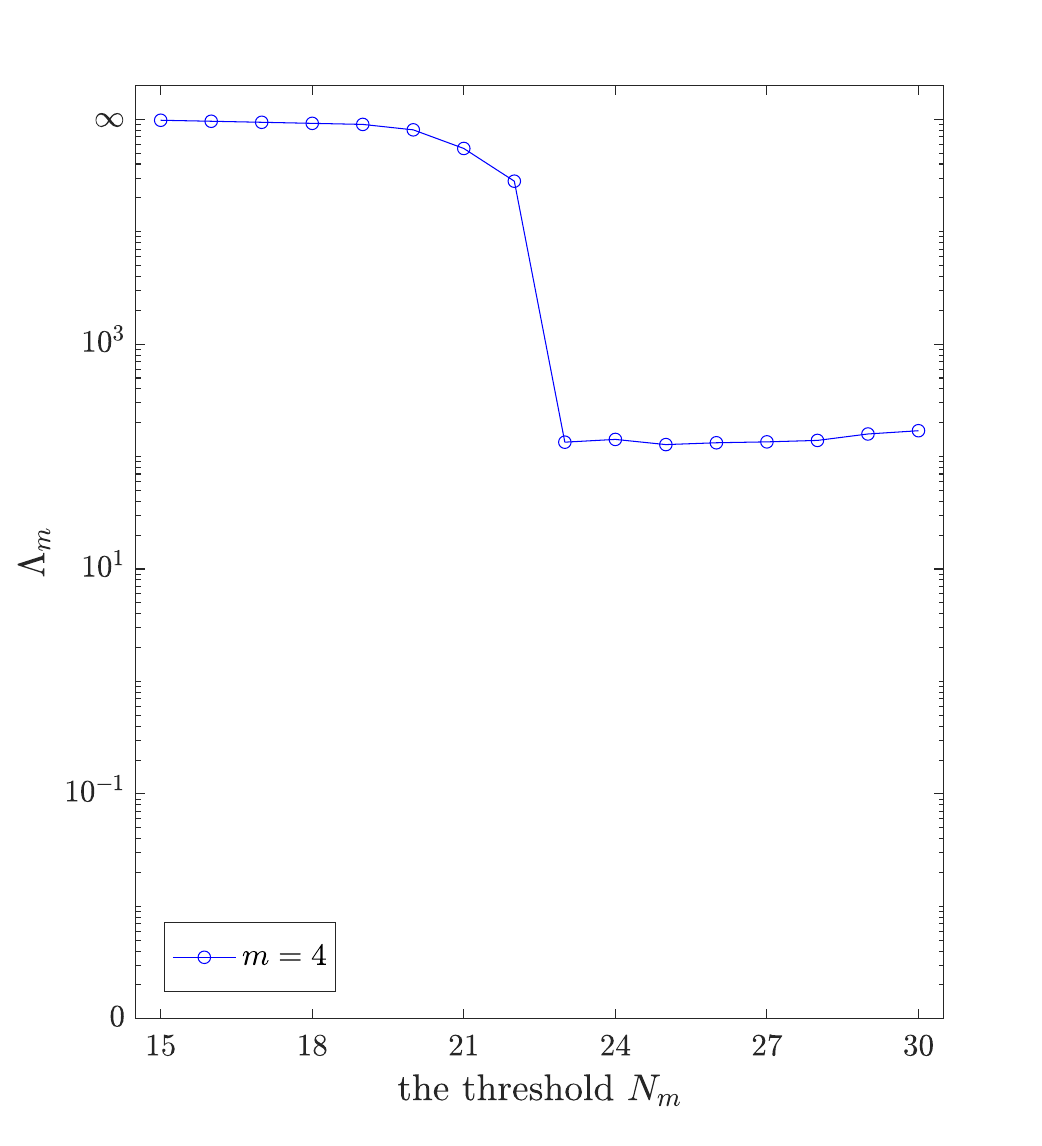}
    \caption{The constant $\Lambda_m$ in two dimensions.}
    \label{fig_Lam2d}
  \end{figure}

  In Fig.~\ref{fig_Lam2d}, we plot the constant $\Lambda_m$ against
  the threshold $N_m$ on the triangular mesh with the mesh size $h =
  1/40$.  
  It is evident that the constant $\Lambda_m$ is nearly a constant and
  increases very slightly when selecting a large threshold $N_m$.
  This numerical observation illustrates the statement given in Remark
  \ref{re_lam} that $\Lambda_m$ will admit a uniform upper bound for
  large $N_m$. 
  Roughly speaking, $N_m$ is close to $1.5\dim(\mP_m(\cdot))$. 
  In the computer implementation, the constant $\Lambda_m$ can also
  serve as an indicator to show whether the element patch is proper. 
  After constructing all element patches, we compute the value of
  $\Lambda_m$. 
  If $\Lambda_m \rightarrow \infty$, we can increase the value of
  $N_m$ until $\Lambda_m \leq 10 \min_{K \in \MTh} \Lambda_{m, K}$.

\begin{figure}[htp]
    \centering
    \includegraphics[width=0.23\textwidth]{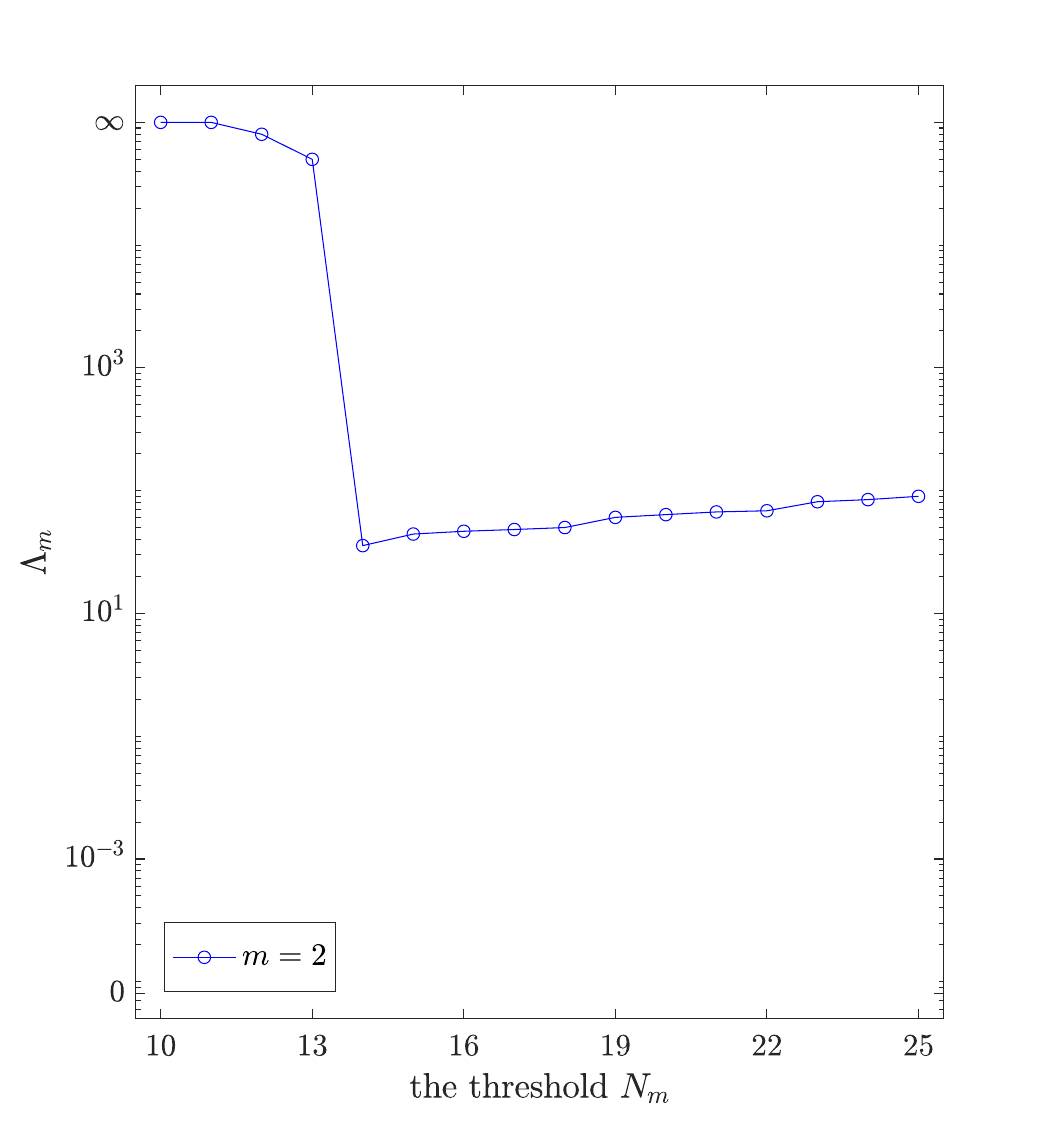}
    \hspace{10pt}
    \includegraphics[width=0.23\textwidth]{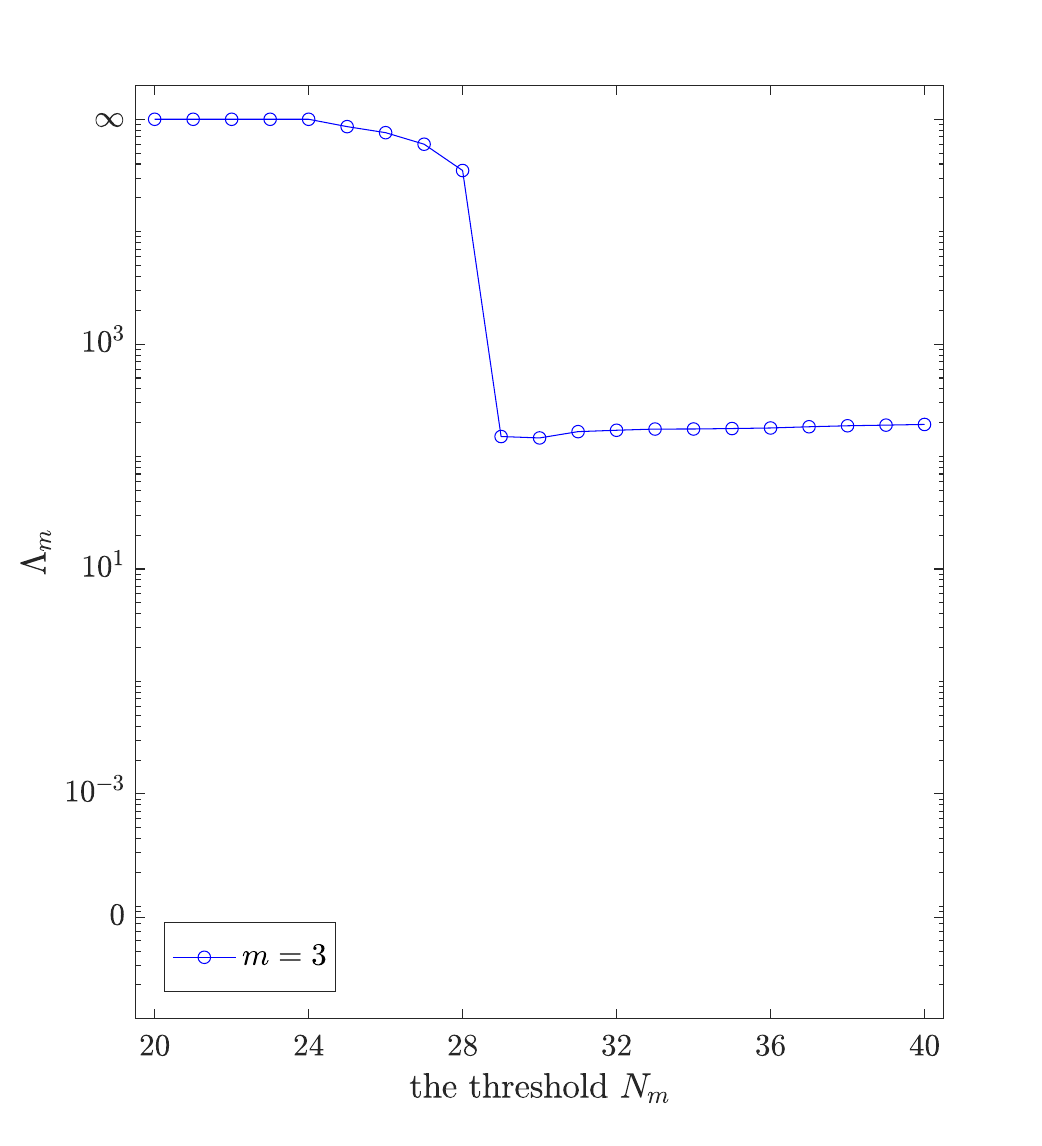}
    \hspace{10pt}
    \includegraphics[width=0.23\textwidth]{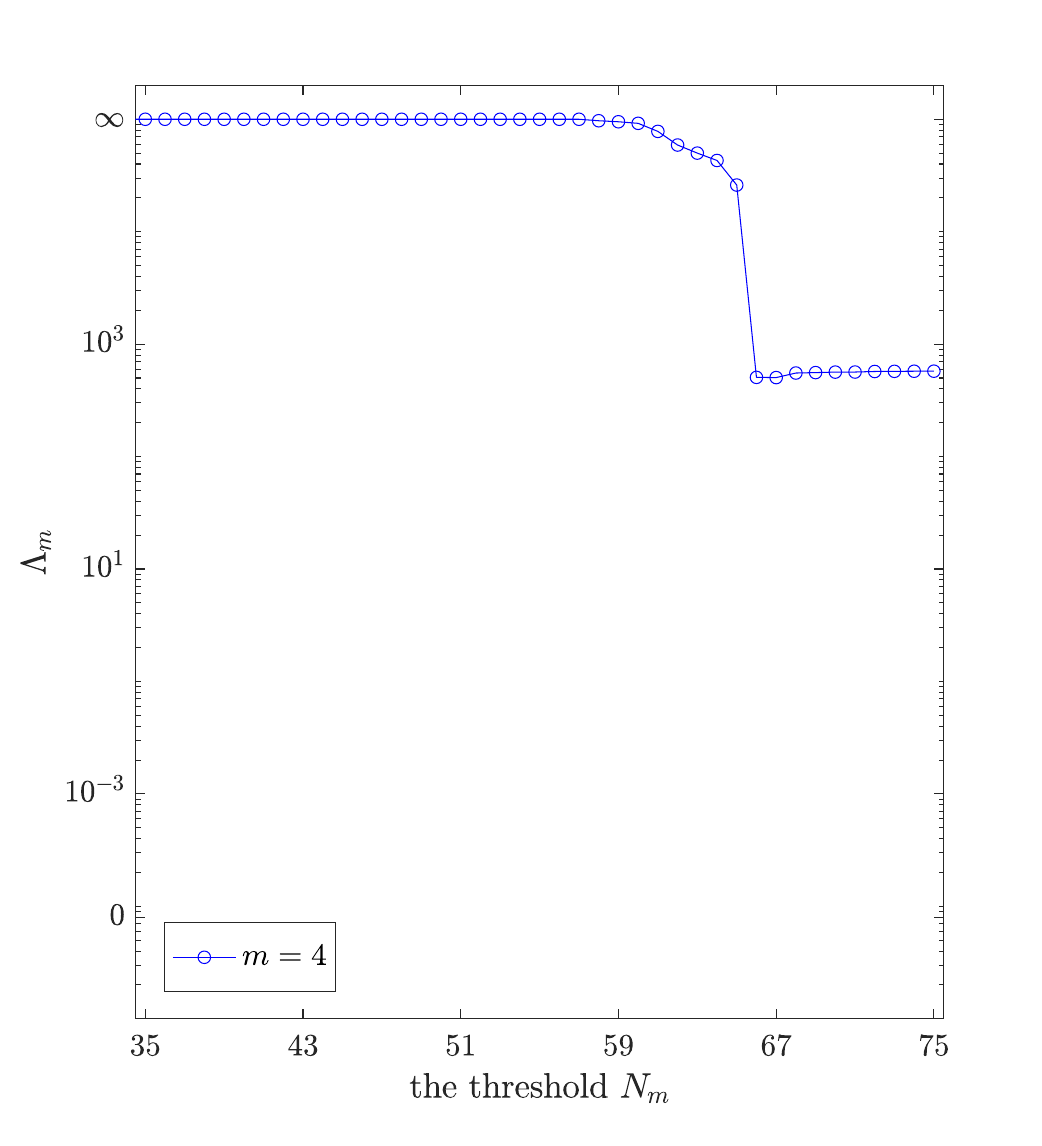}
    \caption{The constant $\Lambda_m$ in three dimensions.}
    \label{fig_Lam3d}
  \end{figure}

  In Fig.~\ref{fig_Lam3d}, we plot the constant $\Lambda_m$ in three
  dimensions on the tetrahedral mesh with the mesh size $h = 1/32$. 
  The numerical result is similar to the case of two dimensions. 
  We still find that $\Lambda_m$ is close to a constant for a large
  $N_m$. 

  \section{}
  \label{sec_app_dMT}
  In this appendix, we verify the estimate \eqref{eq_dMT} in
  both two and three dimensions. 
  Let $W_h^m := \{ w_h \in L^2(\Omega): \ w_h|_K \in \mP_m(K), \
  \forall K \in \MTh \}$ be the discontinuous finite element spaces of
  degree $m$. 
  Let $\bmr{W}_h^m := (W_h^m)^d$, $\bmr{Q}_h^m := \bmr{W}_h^m \cap
  (H^1(\Omega))^d$ be the
  vector-valued discontinuous/$C^0$ finite element spaces. 
  The estimate \eqref{eq_dMT} is based on the following discrete
  Maxwell inequality.
  \begin{lemma}
    There exists a constant $C$ such that
    \begin{equation} 
      \| \nabla \bm{q}_h \|_{L^2(\Omega)}^2 \leq \| \nabla \cdot
      \bm{q}_h \|_{L^2(\Omega)}^2 +  \| \nabla \times \bm{q}_h
      \|_{L^2(\Omega)}^2 + C \sum_{e \in
      \MEhB} h_e^{-1} \| \un \times \bm{q}_h \|_{L^2(e)}^2, \quad
      \forall \bm{q}_h \in \bmr{Q}_h^m.
      \label{eq_dmaxwell}
    \end{equation}
    \label{le_dmaxwell}
  \end{lemma}
  \begin{proof}
    Since $\Omega$ is convex, we first recall the Maxwell inequality
    \cite{Girault1986finite}, which reads 
    \begin{equation}
      \|\nabla {\bm{w}} \|_{L^2(\Omega)} \leq \|\nabla \cdot
      {\bm{w}} \|_{L^2(\Omega)} +  \|\nabla \times {\bm{w}}
      \|_{L^2(\Omega)}, 
      \label{eq_app_maxwellestimate}
    \end{equation}
    for any $\bm{w} \in (H^1(\Omega))^d$ with $\un \times \bm{w} =
    \bm{0}$ on the boundary $\partial \Omega$.
    It is noted that $\bm{q}_h \in (H^1(\Omega))^d$ but its tangential
    trace does not vanish on $\partial \Omega$. 
    We are aiming to construct an auxiliary function $\bm{w}_h$ with
    zero tangential trace on $\partial \Omega$ from $\bm{q}_h$.
    Let $\mc{M} := \{ \bm{\nu}_0, \bm{\nu}_1, \ldots, \bm{\nu}_n\}$
    denote the Lagrange nodes for the space $\bmr{Q}_h^m$, and we let 
    $\{ \phi_{\bm{\nu}_0}, \phi_{\bm{\nu}_1}, \ldots,
    \phi_{\bm{\nu}_n} \}$ be the corresponding basis functions such
    that $\phi_{\bm{\nu}_i}(\bm{\nu}_j) = \delta_{ij}$. 
    We let $\bm{q}_h =  (q_h^1, \ldots, q_h^d)^T$ be expanded by the 
    coefficients $\{ \alpha_{\bm{\nu}}^i \}_{\bm{\nu} \in \mc{M}, 1
    \leq i \leq d}$ such that
    \begin{displaymath}
      \bm{q}_h = \begin{bmatrix}
        q_h^1 \\ \ldots \\ q_h^d 
      \end{bmatrix} = \begin{bmatrix}
        \sum_{\bm{\nu} \in \mc{M}} \alpha_{\bm{\nu}}^1 \phi_{\bm{\nu}}
        \\ 
        \ldots \\
        \sum_{\bm{\nu} \in \mc{M}} \alpha_{\bm{\nu}}^d \phi_{\bm{\nu}} 
      \end{bmatrix} = \sum_{\bm{\nu} \in \mc{M}} \bm{\alpha}_{\bm{\nu}}
      \phi_{\bm{\nu}}, \quad \bm{\alpha}_{\bm{\nu}} =
      (\alpha_{\bm{\nu}}^1, \ldots, \alpha_{\bm{\nu}}^d)^T.
    \end{displaymath}
    The set $\mc{M}$ is further divided into three subsets: 
    \begin{equation}
      \begin{aligned}
        \mc{M}_i &:= \left\{ \bm{\nu} \in \mc{M}: \  \bm{\nu} \text{
        is interior to the domain } \Omega \right\}, \\
        \mc{M}_v &:= \left\{ \begin{aligned}
          & \left\{ \bm{\nu} \in \mc{M}: \  \bm{\nu} \text{ is a
          vertex of the polygonal boundary } \partial \Omega \right\},
          && d = 2, \\
          & \left\{ \bm{\nu} \in \mc{M}: \   \bm{\nu} \text{ lies on
          an edge of the polyhedral boundary } \partial \Omega
          \right\}, && d = 3, \\
        \end{aligned} \right. \\
        \mc{M}_b &:= \mc{M} \backslash ( \mc{M}_i \cup \mc{M}_v). \\
      \end{aligned}
      \label{eq_app_Mdef}
    \end{equation}
    For any $\bm{\nu} \in \mc{M}_b$, there exists a face $e_{\bm{\nu}}
    \in \MEhB$ such that $\bm{\nu} \in \overline{e_{\bm{\nu}}}$. 
    For any $\bm{\nu} \in \mc{M}_v$, there exist two nonparallel faces
    $e_{\bm{\nu}}^1, e_{\bm{\nu}}^2 \in \MEhB$ such that $\bm{\nu} \in
    \overline{e_{\bm{\nu}}^1} \cap \overline{e_{\bm{\nu}}^2}$. 

    We next construct a group of new coefficients $ \{
    \beta_{\bm{\nu}}^i \}_{ \bm{\nu} \in \mc{M}, 1 \leq i \leq d}$ by
    \begin{equation}
      \beta_{\bm{\nu}}^i = \begin{cases}
        \alpha_{\bm{\nu}}^i, & \bm{\nu} \in \mc{M}_i, \\
        \wt{\beta}_{\bm{\nu}}^i, & \bm{\nu} \in \mc{M}_b, \\
        0, & \bm{\nu} \in \mc{M}_v, \\
      \end{cases}
      \quad 1 \leq i \leq d, \quad \bm{\beta}_{\bm{\nu}} =
      (\beta_{\bm{\nu}}^1, \ldots \beta_{\bm{\nu}}^d)^T.
      \label{eq_app_newbeta}
    \end{equation}
    For any $\bm{\nu} \in \mc{M}_b$, $\wt{\beta}_{\bm{\nu}}^i$ is
    determined by the equation 
    \begin{equation}
      \un \times {\bm{\beta}}_{\bm{\nu}} = \bm{0}, \quad
      \un \cdot {\bm{\beta}}_{\bm{\nu}} = \un \cdot
      \bm{\alpha}_{\bm{\nu}},
      \label{eq_app_tildebeta}
    \end{equation}
    where $\un$ is the unit outward normal vector on $e_{\bm{\nu}}$.
    The new polynomial $\bm{w}_h \in \bmr{Q}_h^m$ is constructed by 
    $ \{ \beta_{\bm{\nu}}^i \}_{\bm{\nu} \in \mc{M},1 \leq i \leq d}$
    that 
    \begin{displaymath}
      \bm{w}_h = \begin{bmatrix}
        w_h^1 \\ \ldots \\ w_h^d 
      \end{bmatrix} = \begin{bmatrix}
        \sum_{\bm{\nu} \in \mc{M}} \beta_{\bm{\nu}}^1 \phi_{\bm{\nu}} 
        \\ \ldots \\
        \sum_{\bm{\nu} \in \mc{M}} \beta_{\bm{\nu}}^d \phi_{\bm{\nu}} 
      \end{bmatrix} = \sum_{\bm{\nu} \in \mc{M}} \bm{\beta}_{\bm{\nu}}
      \phi_{\bm{\nu}}.
    \end{displaymath}
    By \eqref{eq_app_Mdef}, we have that 
    \begin{displaymath}
      \|\nabla(\bm{q}_h - {\bm{w}}_h )\|_{L^2(\Omega)}^2 \leq
      C \sum_{\nu \in \mc{M}_b \cup \mc{M}_v} \| \nabla \phi_{\bm{\nu}}
      \|_{L^2(\Omega)}^2 | \bm{\alpha}_{\bm{\nu}} -
      \bm{\beta}_{\bm{\nu}} |^2.
    \end{displaymath}
    From the scaling argument \cite{Karakashian2007convergence}, we
    know that $\|\nabla \phi_{\bm{\nu}} \|_{L^2(\Omega)}^2 \leq C
    h_{e_{\bm{\nu}}}^{d-2}$ for $\forall \bm{\nu} \in \mc{M}_b$. We
    deduce that 
    \begin{align*}
      \sum_{\bm{\nu} \in \mc{M}_b} \| \nabla & \phi_{\bm{\nu}}
      \|_{L^2(\Omega)}^2   | \bm{\alpha}_{\bm{\nu}} -
      \bm{\beta}_{\bm{\nu}} |^2 \leq C  \sum_{\bm{\nu} \in \mc{M}_b}
      h_{e_{\bm{\nu}}}^{d - 2} \left( | \un \times ( \bm{\alpha}_{\bm{\nu}} -
      \bm{\beta}_{\bm{\nu}}) |^2 + |\un \cdot (
      \bm{\alpha}_{\bm{\nu}} - \bm{\beta}_{\bm{\nu}}) |^2
      \right) \\
      & = C  \sum_{\bm{\nu} \in \mc{M}_b} h_{e_{\bm{\nu}}}^{d - 2}  |
      \un \times 
      \bm{\alpha}_{\bm{\nu}} |^2 = C  \sum_{\bm{\nu} \in \mc{M}_b}
      h_{e_{\bm{\nu}}}^{d - 2}  | \un \times \bm{q}_h(\bm{\nu}) |^2
      \\
      & \leq C  \sum_{\bm{\nu} \in \mc{M}_b} h_{e_{\bm{\nu}}}^{d - 2} \|
      \un \times \bm{q}_h \|_{L^\infty(e_{\bm{\nu}})}^2 \leq C
      \sum_{\bm{\nu} \in \mc{M}_b} h_{e_{\bm{\nu}}}^{-1 } \| \un
      \times \bm{q}_h
      \|_{L^2(e_{\bm{\nu}})}^2 \leq  C \sum_{e \in \MEh^b}
      h_e^{-1 } \| \un \times \bm{q}_h \|_{L^2(e)}^2.
    \end{align*}
    For any $\bm{\nu} \in \mc{M}_v$, we let $\un_1, \un_2$ be the unit
    outward normal vectors on $e_{\bm{\nu}}^1, e_{\bm{\nu}}^2$,
    respectively.
    Since $\un_1$ and $\un_2$ are not parallel, by the
    norm equivalence over finite dimensional spaces, there exists a
    constant $C$ only depends on $\Omega$ such that 
    \begin{displaymath}
      |\bm{v}|^2 \leq C ( |\un_1 \times \bm{v}|^2 + |\un_2 \times \bm{v} 
      |^2 ),  \quad \forall \bm{v} \in \mb{R}^d.
    \end{displaymath}
    We derive that 
    \begin{align*}
      \sum_{\bm{\nu} \in \mc{M}_v} &\| \nabla \phi_{\bm{\nu}}
      \|_{L^2(\Omega)}^2  | \bm{\alpha}_{\bm{\nu}} -
      \bm{\beta}_{\bm{\nu}} |^2 =  \sum_{\bm{\nu} \in \mc{M}_v} \|
      \nabla \phi_{\bm{\nu}} \|_{L^2(\Omega)}^2 |
      \bm{\alpha}_{\bm{\nu}}|^2 \leq C  \sum_{\bm{\nu} \in \mc{M}_v}
      \| \nabla \phi_{\bm{\nu}} \|_{L^2(\Omega)}^2   \left(
      |\un_1 \times \bm{\alpha}_{\bm{\nu}}|^2 +
      |\un_2 \times \bm{\alpha}_{\bm{\nu}}|^2 \right) \\ 
      & \leq C \sum_{\bm{\nu} \in \mc{M}_v} \left(
      h_{e_{\bm{\nu}}^1}^{d -2} | \un_1 \times
      \bm{\alpha}_{\bm{\nu}}|^2 + h_{e_{\bm{\nu}}^2}^{d -2} | \un_2
      \times \bm{\alpha}_{\bm{\nu}}|^2 \right)  \leq C \sum_{\bm{\nu}
      \in \mc{M}_v} \left( h_{e_{\bm{\nu}}^1}^{d -2} | \un_1 \times
      \bm{q}_h(\bm{\nu})|^2 + h_{e_{\bm{\nu}}^2}^{d- 2} |\un_2 \times
      \bm{q}_h(\bm{\nu}) |^2 \right) \\
      & \leq C \sum_{\bm{\nu} \in \mc{M}_v} \left(
      h_{e_{\bm{\nu}}^1}^{ -1} \| \un_1 \times  \bm{q}_h \|_{L^2
      (e_{\bm{\nu}}^1)}^2 + h_{e_{\bm{\nu}}^2}^{-1} \| \un_2 \times
      \bm{q}_h \|_{L^2(e_{\bm{\nu}}^2)}^2 \right) \leq  C \sum_{e \in
      \MEh^b} h_e^{-1 } \| \un \times \bm{q}_h \|_{L^2(e)}^2.   
    \end{align*}
    Finally, we arrive at
    \begin{equation}
      \| \nabla (\bm{q}_h - \bm{w}_h) \|_{L^2(\Omega)}^2 \leq C
      \sum_{e \in \MEhB} h_e^{-1} \| \un \times \bm{q}_h
      \|_{L^2(e)}^2,
      \label{eq_qtqdiff}
    \end{equation}
    Clearly, $\bm{w}_h$ satisfies the estimate
    \eqref{eq_app_maxwellestimate}. 
    Together with \eqref{eq_qtqdiff}, we have that 
    \begin{align*}
      \| \nabla \bm{q}_h \|_{L^2(\Omega)} & \leq \| \nabla \bm{w}_h
      \|_{L^2(\Omega)} + \| \nabla(\bm{q}_h -
      \bm{w}_h)\|_{L^2(\Omega)} \\
      & \leq \|\nabla \cdot {\bm{w}}_h \|_{L^2(\Omega)} +  \|\nabla
      \times {\bm{w}}_h \|_{L^2(\Omega)} +  \| \nabla(\bm{q}_h -
      \bm{w}_h)\|_{L^2(\Omega)} \\ 
      & \leq \| \nabla \cdot \bm{q}_h \|_{L^2(\Omega)} + \| \nabla
      \times \bm{q}_h \|_{L^2(\Omega)} + C  \| \nabla(\bm{q}_h -
      \bm{w}_h)\|_{L^2(\Omega)},
    \end{align*}
    which indicates the estimate \eqref{eq_dmaxwell}. This completes
    the proof.
  \end{proof}
  The estimate \eqref{eq_dmaxwell} can be extended to the
  discontinuous space $\bmr{W}_h^m$, which reads 
  \begin{equation}
    \begin{aligned}
      \sum_{K \in \MTh} \| \nabla \bm{w}_h \|_{L^2(K)}^2 & \leq
      \sum_{K \in \MTh} ( \| \nabla \cdot \bm{w}_h \|_{L^2(K)}^2 + \|
      \nabla \times \bm{w}_h \|_{L^2(K)}^2) \\
      + & C \sum_{e \in \MEhI} h_e^{-1} \| \jump{ \bm{w}_h}
      \|_{L^2(e)}^2 + C \sum_{e \in \MEhB} h_e^{-1} \| \un \times
      \bm{w}_h \|_{L^2(e)}^2, \quad \forall \bm{w}_h \in \bmr{W}_h^m,
    \end{aligned}
    \label{eq_dmaxwellwh}
  \end{equation}
  by applying the Oswald interpolant. 
  From \cite[Theorem 2.1]{Karakashian2007convergence}, for any
  $\bm{w}_h \in \bmr{W}_h^m$, there exists $\bm{q}_h \in \bmr{Q}_h^m$
  such that
  \begin{equation}
    \sum_{K \in \MTh} ( h_K^2 \| \bm{w}_h - \bm{q}_h \|_{L^2(K)}^2 + \|
    \bm{w}_h - \bm{q}_h \|_{H^1(K)}^2) \leq C \sum_{e \in \MEhI}
    h_e^{-1} \| \jump{\bm{w}_h} \|_{L^2(e)}^2.
    \label{eq_whdiffqh}
  \end{equation}
  Then, the estimate \eqref{eq_dmaxwellwh} can be easily verified
  using estimates \eqref{eq_whdiffqh} and \eqref{eq_dmaxwell} and the
  trace estimate.


  Finally, let us present the main conclusion in this appendix. 
  \begin{lemma}
    There exists a constant $C$ such that
    \begin{align}
      \sum_{K \in \MTh} \| D^2 v_h \|_{L^2(K)}^2 & \leq 
      \sum_{K \in \MTh} \| \Delta v_h \|_{L^2(K)}^2 
      + C ( \sum_{e \in \MEhI} h_e^{-1} \| \jump{\nabla_{\un} v_h}
      \|_{L^2(e)}^2 + \sum_{e \in \MEh} h_e^{-3} \| \jump{v_h}
      \|_{L^2(e)}^2), 
      \label{eq_D2vhbound}
    \end{align}
    for $\forall v_h \in W_h^m$.
    \label{le_D2vhbound}
  \end{lemma}
  \begin{proof}
    Let $\bm{q}_h \in \bmr{W}_h^m$ be determined by $\bm{q}_h|_K =
    \nabla v_h|_K$ for $\forall K \in \MTh$.
    Then, it is clear that $\nabla \times \bm{q}_h|_K = 0$ on every
    $K \in \MTh$.
    For any interior face $e \in \MEhI$, we apply the
    inverse estimate to find that
    \begin{align*}
      h_e^{-1} \| \jump{\bm{q}_h} \|_{L^2(e)}^2 & \leq h_e^{-1} ( \|
      \jump{\un \times \bm{q}_h} \|_{L^2(e)}^2 +  \|
      \jump{\un \cdot \bm{q}_h} \|_{L^2(e)}^2) =  h_e^{-1} ( \|
      \jump{\un \times \nabla v_h} \|_{L^2(e)}^2 +  \|
      \jump{\nabla_{\un} v_h} \|_{L^2(e)}^2) \\ 
      & \leq C( h_e^{-3} \| \jump{v_h} \|_{L^2(e)}^2 + h_e^{-1} \|
      \jump{\nabla_{\un} v_h} \|_{L^2(e)}^2).
    \end{align*}
    From the estimate \eqref{eq_dmaxwellwh}, we have that
    \begin{align*}
      \sum_{K \in \MTh} &\| D^2 v_h \|_{L^2(K)}^2  = \sum_{K \in \MTh}
      \| \nabla \bm{q}_h \|_{L^2(K)}^2 \\
      & \leq \sum_{K \in \MTh} \|
      \nabla \cdot \bm{q}_h  \|_{L^2(K)}^2 + C \sum_{e \in \MEhI}
      h_e^{-1}  \| \jump{\bm{q}_h} \|_{L^2(e)}^2  + C \sum_{e
      \in \MEhB} h_e^{-1} \| \un \times \bm{q}_h \|_{L^2(e)}^2 \\
      & \leq \sum_{K \in \MTh} \|
      \Delta v_h  \|_{L^2(K)}^2  + C \sum_{e \in \MEhI}  h_e^{-1}  \|
      \jump{\nabla v_h} \|_{L^2(e)}^2 + C \sum_{e \in \MEh} h_e^{-3}
      \| \jump{v_h} \|_{L^2(e)}^2,
    \end{align*}
    which completes the proof.
  \end{proof}
\end{appendix}

\bibliographystyle{amsplain}
\bibliography{../ref}

\providecommand{\bysame}{\leavevmode\hbox to3em{\hrulefill}\thinspace}
\providecommand{\MR}{\relax\ifhmode\unskip\space\fi MR }
\providecommand{\MRhref}[2]{%
  \href{http://www.ams.org/mathscinet-getitem?mr=#1}{#2}
}
\providecommand{\href}[2]{#2}
\begin{thebibliography}{10}

\bibitem{Adini1960analysis}
A.~Adini and R.~W Clough, \emph{{Analysis of Plate Bending by the Finite
  Element Method}}, NSF report, 1961.

\bibitem{Argyris1968tuba}
J.~H. Argyris, I.~Fried, and D.~W Scharpf, \emph{The {TUBA} family of plate
  elements for the matrix displacement method}, Aeronaut. J., Roy Aeronaut.
  Soc. \textbf{72} (1968), no.~692, 701--709.

\bibitem{Blum1980boundary}
H.~Blum and R.~Rannacher, \emph{On the boundary value problem of the biharmonic
  operator on domains with angular corners}, Math. Methods Appl. Sci.
  \textbf{2} (1980), no.~4, 556--581.

\bibitem{Boffi2013mixed}
D.~Boffi, F.~Brezzi, and M.~Fortin, \emph{{Mixed Finite Element Methods and
  Applications}}, Springer Series in Computational Mathematics, vol.~44,
  Springer, Heidelberg, 2013.

\bibitem{Brenner1989optimal}
S.~C. Brenner, \emph{An optimal-order nonconforming multigrid method for the
  biharmonic equation}, SIAM J. Numer. Anal. \textbf{26} (1989), no.~5,
  1124--1138.

\bibitem{Brenner1996two}
\bysame, \emph{Two-level additive {S}chwarz preconditioners for nonconforming
  finite element methods}, Math. Comp. \textbf{65} (1996), no.~215, 897--921.

\bibitem{Brenner2015C0}
S.~C. Brenner, P.~Monk, and J.~Sun, \emph{{$C^0$} interior penalty {G}alerkin
  method for biharmonic eigenvalue problems}, Spectral and High Order Methods
  for Partial Differential Equations---{ICOSAHOM} 2014, Lect. Notes Comput.
  Sci. Eng., vol. 106, Springer, Cham, 2015, pp.~3--15.

\bibitem{Brenner2005interior}
S.~C. Brenner and L.-Y. Sung, \emph{{$C^0$} interior penalty methods for fourth
  order elliptic boundary value problems on polygonal domains}, J. Sci. Comput.
  \textbf{22/23} (2005), 83--118.

\bibitem{Brenner2006multigrid}
\bysame, \emph{Multigrid algorithms for {$C^0$} interior penalty methods}, SIAM
  J. Numer. Anal. \textbf{44} (2006), no.~1, 199--223.

\bibitem{Brenner2005two}
S.~C. Brenner and K.~Wang, \emph{Two-level additive {S}chwarz preconditioners
  for {$C^0$} interior penalty methods}, Numer. Math. \textbf{102} (2005),
  no.~2, 231--255.

\bibitem{Carstensen2021hierarchical}
C.~Carstensen and J.~Hu, \emph{Hierarchical {A}rgyris finite element method for
  adaptive and multigrid algorithms}, Comput. Methods Appl. Math. \textbf{21}
  (2021), no.~3, 529--556.

\bibitem{Chalmers2018low}
N.~Chalmers and T.~Warburton, \emph{Low-order preconditioning of high-order
  triangular finite elements}, SIAM J. Sci. Comput. \textbf{40} (2018), no.~6,
  A4040--A4059.

\bibitem{Chen2017C0}
H.~Chen, H.~Guo, Z.~Zhang, and Q.~Zou, \emph{A {$C^0$} linear finite element
  method for two fourth-order eigenvalue problems}, IMA J. Numer. Anal.
  \textbf{37} (2017), no.~4, 2120--2138.

\bibitem{Ciarlet2002finite}
P.~G. Ciarlet, \emph{{The Finite Element Method for Elliptic Problems}},
  Classics in Applied Mathematics, vol.~40, Society for Industrial and Applied
  Mathematics (SIAM), Philadelphia, PA, 2002, Reprint of the 1978 original
  [North-Holland, Amsterdam; MR0520174 (58 \#25001)].

\bibitem{Cockburn2009hybridizable}
B.~Cockburn, B.~Dong, and J.~Guzm\'an, \emph{A hybridizable and superconvergent
  discontinuous {G}alerkin method for biharmonic problems}, J. Sci. Comput.
  \textbf{40} (2009), no.~1-3, 141--187.

\bibitem{Douglas1979family}
J.~Douglas, Jr., T.~Dupont, P.~Percell, and R.~Scott, \emph{A family of
  {$C\sp{1}$} finite elements with optimal approximation properties for various
  {G}alerkin methods for 2nd and 4th order problems}, RAIRO Anal. Num\'er.
  \textbf{13} (1979), no.~3, 227--255.

\bibitem{Engel2002continuous}
G.~Engel, K.~Garikipati, T.~J.~R. Hughes, M.~G. Larson, L.~Mazzei, and R.~L.
  Taylor, \emph{Continuous/discontinuous finite element approximations of
  fourth-order elliptic problems in structural and continuum mechanics with
  applications to thin beams and plates, and strain gradient elasticity},
  Comput. Methods Appl. Mech. Engrg. \textbf{191} (2002), no.~34, 3669--3750.

\bibitem{Feng2005two}
X.~Feng and O.~A. Karakashian, \emph{Two-level non-overlapping {S}chwarz
  preconditioners for a discontinuous {G}alerkin approximation of the
  biharmonic equation}, J. Sci. Comput. \textbf{22/23} (2005), 289--314.

\bibitem{Georgoulis2008discontinuous}
E.~H. Georgoulis and P.~Houston, \emph{Discontinuous {G}alerkin methods for the
  biharmonic problem}, IMA J. Numer. Anal. \textbf{29} (2009), no.~3, 573--594.

\bibitem{Georgoulis2011posteriori}
E.~H. Georgoulis, P.~Houston, and J.~Virtanen, \emph{An {\it a posteriori}
  error indicator for discontinuous {G}alerkin approximations of fourth-order
  elliptic problems}, IMA J. Numer. Anal. \textbf{31} (2011), no.~1, 281--298.

\bibitem{Girault1986finite}
V.~Girault and P.-A. Raviart, \emph{Finite element methods for
  {N}avier-{S}tokes equations}, Springer Series in Computational Mathematics,
  vol.~5, Springer-Verlag, Berlin, 1986, Theory and algorithms.

\bibitem{Gudi2008mixed}
T.~Gudi, N.~Nataraj, and A.~K. Pani, \emph{Mixed discontinuous {G}alerkin
  finite element method for the biharmonic equation}, J. Sci. Comput.
  \textbf{37} (2008), no.~2, 139--161.

\bibitem{Guo2018C0}
H.~Guo, Z.~Zhang, and Q.~Zou, \emph{A {$C^0$} linear finite element method for
  biharmonic problems}, J. Sci. Comput. \textbf{74} (2018), no.~3, 1397--1422.

\bibitem{Huang2020recovery}
Y.~Huang, H.~Wei, W.~Yang, and N.~Yi, \emph{Recovery based finite element
  method for biharmonic equation in 2{D}}, J. Comput. Math. \textbf{38} (2020),
  no.~1, 84--102.

\bibitem{Karakashian2018two}
O.~A. Karakashian and C.~Collins, \emph{Two-level additive {S}chwarz methods
  for discontinuous {G}alerkin approximations of the biharmonic equation}, J.
  Sci. Comput. \textbf{74} (2018), no.~1, 573--604.

\bibitem{Karakashian2007convergence}
O.~A. Karakashian and F.~Pascal, \emph{Convergence of adaptive discontinuous
  {G}alerkin approximations of second-order elliptic problems}, SIAM J. Numer.
  Anal. \textbf{45} (2007), no.~2, 641--665.

\bibitem{Lamichhane2014finite}
B.~P. Lamichhane, \emph{A finite element method for a biharmonic equation based
  on gradient recovery operators}, BIT \textbf{54} (2014), no.~2, 469--484.

\bibitem{Li2023preconditioned}
R.~Li, Q.~Liu, and F.~Yang, \emph{Preconditioned nonsymmetric/symmetric
  discontinuous {G}alerkin method for elliptic problem with reconstructed
  discontinuous approximation}, accpeted by J. Sci. Comput. (2023).

\bibitem{Li2017discontinuous}
R.~Li, P.~Ming, Z.~Sun, F.~Yang, and Z.~Yang, \emph{A discontinuous {G}alerkin
  method by patch reconstruction for biharmonic problem}, J. Comput. Math.
  \textbf{37} (2019), no.~4, 563--580.

\bibitem{Li2016discontinuous}
R.~Li, P.~Ming, Z.~Sun, and Z.~Yang, \emph{An arbitrary-order discontinuous
  {G}alerkin method with one unknown per element}, J. Sci. Comput. \textbf{80}
  (2019), no.~1, 268--288.

\bibitem{Li2012efficient}
R.~Li, P.~Ming, and F.~Tang, \emph{An efficient high order heterogeneous
  multiscale method for elliptic problems}, Multiscale Model. Simul.
  \textbf{10} (2012), no.~1, 259--283.

\bibitem{Morley1968triangular}
L.~Morley, \emph{The triangular equilibrium element in the solution of plate
  bending problems}, Aero. Quart. \textbf{19} (1968), no.~2, 149--169.

\bibitem{Mozolevski2003priori}
I.~Mozolevski and E.~S\"uli, \emph{A priori error analysis for the
  {$hp$}-version of the discontinuous {G}alerkin finite element method for the
  biharmonic equation}, Comput. Methods Appl. Math. \textbf{3} (2003), no.~4,
  596--607.

\bibitem{Mozolevski2007hp}
I.~Mozolevski, E.~S\"uli, and P.~R. B\"osing, \emph{{$hp$}-version a priori
  error analysis of interior penalty discontinuous {G}alerkin finite element
  approximations to the biharmonic equation}, J. Sci. Comput. \textbf{30}
  (2007), no.~3, 465--491.

\bibitem{Neilan2019discrete}
M.~Neilan and M.~Wu, \emph{Discrete {M}iranda-{T}alenti estimates and
  applications to linear and nonlinear {PDE}s}, J. Comput. Appl. Math.
  \textbf{356} (2019), 358--376.

\bibitem{Onate1993derivation}
E.~O\~nate and M.~Cervera, \emph{Derivation of thin plate bending elements with
  one degree of freedom per node: a simple three node triangle}, Engrg. Comput.
  \textbf{10} (1993), no.~6, 543--561.

\bibitem{Pazner2023low}
W.~Pazner, T.~Kolev, and C.~R. Dohrmann, \emph{Low-order preconditioning for
  the high-order finite element de {R}ham complex}, SIAM J. Sci. Comput.
  \textbf{45} (2023), no.~2, A675--A702.

\bibitem{Phaal1992simple}
R.~Phaal and C.~R. Calladine, \emph{A simple class of finite elements for plate
  and shell problems. {II}: An element for thin shells, with only translational
  degrees of freedom}, Internat. J. Numer. Methods Engrg. \textbf{35} (1992),
  no.~5, 979--996.

\bibitem{Powell1981approximation}
M.~J.~D. Powell, \emph{Approximation theory and methods}, Cambridge University
  Press, Cambridge-New York, 1981.

\bibitem{Rannacher1979nonconforming}
R.~Rannacher, \emph{On nonconforming a mixed finite element methods for plate
  bending problems. {T}he linear case}, RAIRO. Anal. numér. \textbf{13}
  (1979), no.~4, 369--387.

\bibitem{Shi1986convergence}
Z.~C. Shi, \emph{On the convergence of the incomplete biquadratic nonconforming
  plate element}, Math. Numer. Sinica \textbf{8} (1986), no.~1, 53--62.
  \MR{864031}

\bibitem{Smears2018nonoverlapping}
I.~Smears, \emph{Nonoverlapping domain decomposition preconditioners for
  discontinuous {G}alerkin approximations of {H}amilton-{J}acobi-{B}ellman
  equations}, J. Sci. Comput. \textbf{74} (2018), no.~1, 145--174.

\bibitem{Smears2014discontinuous}
I.~Smears and E.~S\"{u}li, \emph{Discontinuous {G}alerkin finite element
  approximation of {H}amilton-{J}acobi-{B}ellman equations with {C}ordes
  coefficients}, SIAM J. Numer. Anal. \textbf{52} (2014), no.~2, 993--1016.

\bibitem{Stevenson2003analysis}
R.~Stevenson, \emph{An analysis of nonconforming multi-grid methods, leading to
  an improved method for the {M}orley element}, Math. Comp. \textbf{72} (2003),
  no.~241, 55--81.

\bibitem{Suli2007hp}
E.~S\"uli and I.~Mozolevski, \emph{{$hp$}-version interior penalty {DGFEM}s for
  the biharmonic equation}, Comput. Methods Appl. Mech. Engrg. \textbf{196}
  (2007), no.~13-16, 1851--1863.

\bibitem{Tang2017local}
S.~Tang and X.~Xu, \emph{Local multilevel methods with rectangular finite
  elements for the biharmonic problem}, SIAM J. Sci. Comput. \textbf{39}
  (2017), no.~6, A2592--A2615.

\bibitem{Vanek2001convergence}
P.~Van\v{e}k, M.~Brezina, and J.~Mandel, \emph{Convergence of algebraic
  multigrid based on smoothed aggregation}, Numer. Math. \textbf{88} (2001),
  no.~3, 559--579.

\bibitem{Xu1992iterative}
J.~Xu, \emph{Iterative methods by space decomposition and subspace correction},
  SIAM Rev. \textbf{34} (1992), no.~4, 581--613.

\bibitem{Zhang1989optimal}
S.~Zhang, \emph{An optimal order multigrid method for biharmonic, {$C^1$}
  finite element equations}, Numer. Math. \textbf{56} (1989), no.~6, 613--624.

\bibitem{Zhang1996two}
X.~Zhang, \emph{Two-level {S}chwarz methods for the biharmonic problem
  discretized conforming {$C^1$} elements}, SIAM J. Numer. Anal. \textbf{33}
  (1996), no.~2, 555--570.

\end{thebibliography}

\end{document}